\documentclass[11pt]{amsart}
\usepackage[utf8]{inputenc}
\usepackage{amssymb}
\usepackage{amsmath}
\usepackage{mathtools}
\usepackage{ifpdf}
\usepackage{vmargin}

\usepackage{graphicx}  
\usepackage{subfigure}
\usepackage{epstopdf}
\usepackage{caption}

\usepackage{setspace}
\usepackage{xcolor}

\ifpdf 
\usepackage[pdftex]{hyperref}
\else
\usepackage[hypertex]{hyperref}
\fi

\hypersetup{
	hidelinks,
	pdftitle={Sharp global well-posedness for the defocusing nonlinear wave equation outside a ball},
	pdfauthor={Guixiang Xu},
}

\newtheorem{theorem}{Theorem}
\newtheorem{proposition}[theorem]{Proposition}
\newtheorem{lemma}[theorem]{Lemma}
\newtheorem{remark}[theorem]{Remark}
\newtheorem{corollary}[theorem]{Corollary}

\numberwithin{equation}{section}
\newcommand{\R}{\mathbb{R}}

\newcommand{\bee}{\begin{eqnarray*}}
	\newcommand{\eee}{\end{eqnarray*}}

\numberwithin{equation}{section}
\numberwithin{theorem}{section}

\begin{document}

% \doublespacing
\onehalfspacing

\title[Sharp global well-posedness for the defocusing NLW outside a ball]{Sharp global well-posedness for the defocusing nonlinear wave equation outside a ball with radial data}

\date\today

\author[Guixiang Xu]{Guixiang Xu}
\address{\hskip-1.15em Guixiang Xu 
	\hfill\newline Laboratory of Mathematics and Complex Systems,
	\hfill\newline Ministry of Education,
	\hfill\newline School of Mathematical Sciences,
	\hfill\newline Beijing Normal University,
	\hfill\newline Beijing, 100875, People's Republic of China.}
\email{guixiang@bnu.edu.cn}


\begin{abstract}
Let $\Omega=\mathbb{R}^3\setminus\overline{B(0,1)}$ and $3<p<5$. We prove sharp global well-posedness of  the
defocusing Dirichlet boundary value problem
of nonlinear wave equation with radial data in the critical space
$\dot H^{s_c}_{D,\mathrm{rad}}(\Omega)\times
\dot H^{s_c-1}_{D,\mathrm{rad}}(\Omega)$, where
$s_c=\frac32-\frac{2}{p-1}$, thereby improve the result in
\cite{XuYangYou2025}. The proof combines the distorted Fourier transform in
\cite{LiSZ:NLS,XuYang2024} with the energy--Morawetz estimate  by Dodson in 
\cite{DodsonCubic2018,Dodson2018}. 
\end{abstract}

\maketitle

\section{Introduction}\label{sec:introduction}

We consider
\begin{equation}
 \left\{
 \begin{aligned}
  &\partial_t^2u-\Delta u+F(u)=0,
       &&(t,x)\in I\times\Omega_R,\\
  &u(t,x)=0,&&x\in\partial\Omega_R,\\
  &(u,\partial_tu)|_{t=0}=(u_0,u_1),
 \end{aligned}
 \right.
 \qquad
 \label{eq:nlw}
\end{equation}
where $F(u):=\left\lvert u\right\rvert^{p-1}u$, $3<p<5$, the solution $u$ is radial and real-valued function, and
\begin{equation}
\Omega_R:=\{x\in\mathbb{R}^3:\left\lvert x\right\rvert>R\},
\qquad R>0.
\label{eq:domain}
\end{equation}
We consider the scale-covariant family $\{\Omega_R\}_{R>0}$ and return to
the unit exterior domain in Theorem \ref{thm:main} and at the end of the paper. The functional calculus gives
\begin{equation*}
u(t)=\mathcal W_R(t)(u_0,u_1)
-\int_0^t\frac{\sin((t-s)\sqrt{-\Delta_{\Omega_R}})}
{\sqrt{-\Delta_{\Omega_R}}}F(u(s))\,\mathrm{d}s,
\end{equation*}
where
$
\mathcal W_R(t)(u_0,u_1)
=\cos(t\sqrt{-\Delta_{\Omega_R}})u_0
+\frac{\sin(t\sqrt{-\Delta_{\Omega_R}})}
{\sqrt{-\Delta_{\Omega_R}}}u_1.
$
For sufficiently regular solutions, the conserved energy is
\begin{align*}
E(u)(t):=& \int_{\Omega_R}\left(
\frac12|\nabla u(t,x)|^2+\frac12|\partial_tu(t,x)|^2
+\frac1{p+1}|u(t,x)|^{p+1}\right)\,\mathrm{d}x,\\
 =& \, E(u)(0),\qquad t\in I.
\end{align*}

Although a fixed exterior domain is not invariant under scaling, the family
in \eqref{eq:domain} is scale covariant. More precisely, if $u$ solves
\eqref{eq:nlw} on $\Omega_1$, then
\begin{equation}
 u^{[\lambda]}(t,x)
 :=\lambda^{\frac{2}{p-1}}u(\lambda t,\lambda x)
 \label{eq:scaling}
\end{equation}
solves the same equation on $\Omega_{1/\lambda}$. The critical norm
$\dot H_D^{s_c}\times\dot H_D^{s_c-1}$ with $s_c=\frac32-\frac{2}{p-1}$ is invariant under
\eqref{eq:scaling}. Indeed, if
$f^{[\lambda]}(x)=\lambda^{2/(p-1)}f(\lambda x)$, then
\begin{equation*}
 \left\lVert f^{[\lambda]}\right\rVert_{\dot H^s}
 =\lambda^{s-\frac32+\frac{2}{p-1}}\left\lVert f\right\rVert_{\dot H^s}.
\end{equation*}
The exponent vanishes if $s=s_c$, which is the reason why we call \eqref{eq:nlw} $\dot H^{s_c}$-critical. By contrast, the
energy functional satisfies
\begin{equation*}
 E\bigl(u^{[\lambda]}(0)\bigr)
 =\lambda^{\frac{5-p}{p-1}}E(u(0)).
\end{equation*}
Thus the equation is energy critical for $p=5$, while for $3<p<5$, the equation is energy
subcritical. The main difficulty is the long-time propagation and
control of solution of \eqref{eq:nlw} with low-regularity initial data which may have infinite
(kinetic) energy.

The corresponding Cauchy problem with radial data on the whole space is now well understood in the critical
Sobolev space. For the defocusing, cubic NLW equation, Dodson proved in
\cite{DodsonCubic2018} that every solution with radial data
\begin{equation*}
 (u_0,u_1)\in
 \dot H^{1/2}(\mathbb{R}^3)\times\dot H^{-1/2}(\mathbb{R}^3)
\end{equation*}
exists globally and scatters in both time directions.
Moreover, its critical spacetime norm obeys a qualitative bound
\begin{equation*}
 \left\|u\right\|_{L^4_{t,x}(\mathbb{R}\times\mathbb{R}^3)}
 \leq C\!\left(
 \left\|u_0\right\|_{\dot H^{1/2}(\mathbb{R}^3)}
 +\left\|u_1\right\|_{\dot H^{-1/2}(\mathbb{R}^3)}\right).
\end{equation*}
For every $3<p<5$, the same author subsequently established in
\cite{Dodson2018} global well-posedness and scattering for all radial data
in
\begin{equation*}
 \dot H^{s_c}(\mathbb{R}^3)\times
 \dot H^{s_c-1}(\mathbb{R}^3),
 \qquad s_c=\frac32-\frac{2}{p-1},
\end{equation*}
together with the critical spacetime estimate
\begin{equation*}
 \left\|u\right\|_{L^{2(p-1)}_{t,x}
 (\mathbb{R}\times\mathbb{R}^3)}
 \leq C_p\!\left(
 \left\|u_0\right\|_{\dot H^{s_c}(\mathbb{R}^3)}
 +\left\|u_1\right\|_{\dot H^{s_c-1}(\mathbb{R}^3)}\right).
\end{equation*}
Here and below, the functions denoted by $C$ and $C_p$ need not be
explicit. For the related higher-dimensional results, see
\cite{Dodson2024,Dodson2025}.  We can refer to \cite{BahChe:NLW, DodsonCubic2019, GallPlan03:NLW, KenigPV00:NLW,Roy:NLW:radial,Roy:NLW} in the subcritical space. 

These whole-space results are fundamental and essential, but
their proofs cannot be transferred directly to an exterior domain. The
Fourier truncation method in \cite{Bourgain1998} and hyperbolic-coordinate transform in
\cite{DodsonCubic2018,Dodson2018} use translation-invariant Fourier
multipliers and the geometry structure of the light cone. A Dirichlet boundary
replaces convolution kernels of the Littlewood-Paley projection by direct-minus-reflected kernels in \eqref{eq:sinekernel}, and produces a
boundary flux in the Morawetz multiplier identities \eqref{eq:Morawetzidentity}, and breaks the fixed-domain scaling.

The purpose of this paper is to recover sharp  global
well-posedness in $\mathbb{R}^3$ in \cite{Dodson2018} for the exterior Dirichlet problem. The low-regularity global well-posedness outside a ball was initiated in
\cite{XuYang2024}. For the cubic equation, global exotic Strichartz estimates with radial data
combined with a Fourier truncation method in the spirit of
\cite{Bourgain1998} yield global well-posedness for
$3/4<s<1$, subject to an additional $L^4$ condition. The powers $3<p<5$
were subsequently treated in \cite{XuYangYou2025} in the range
\begin{equation*}
 1-\frac{(p+3)(1-s_c)}{4(2p-3)}<s<1.
\end{equation*}

Although there is no  conserved quantity at the critical regularity, the origin point is only possible singular point under the radial assumption; compare
with \cite{Bourgain1999,LiSZ:NLS}. While the exterior domain removes the origin
point. Based on the above heuristic,  the following sharp global well-posedness of \eqref{eq:nlw} in the critical space is expected.

\begin{theorem}\label{thm:main}
Let $3<p<5$. For every pair of real-valued radial initial data
\begin{equation}
 (u_0,u_1)\in
 \dot H^{s_c}_{D,\mathrm{rad}}(\Omega_1)
 \times
 \dot H^{s_c-1}_{D,\mathrm{rad}}(\Omega_1),
 \label{eq:maindata}
\end{equation}
the equation \eqref{eq:nlw} has a unique global solution satisfying
\begin{equation}
 (u,\partial_tu)\in
 C\!\left(\mathbb{R};
 \dot H^{s_c}_{D,\mathrm{rad}}(\Omega_1)
 \times
 \dot H^{s_c-1}_{D,\mathrm{rad}}(\Omega_1)\right).
 \label{eq:mainconclusion}
\end{equation}
For every $T>0$, the data-to-solution map is locally Lipschitz from the
space in \eqref{eq:maindata} into the corresponding continuous solution
space on $[-T,T]$.
\end{theorem}

At the critical level, no additional potential energy assumption is needed.
Indeed, the radial Sobolev estimate in
Lemma~\ref{lem:radialsobolev} gives
\begin{equation}
\dot H^{s_c}_{D,\mathrm{rad}}(\Omega_1)
\hookrightarrow L^{p+1}(\Omega_1),
\qquad 3<p<5.
\label{eq:criticalLp}
\end{equation}
The solution in Theorem~\ref{thm:main} may have infinite
energy because of $s_c<1$. The cubic endpoint case $s=1/2$ was obtained more recently in \cite{Lin2026} by
means of the radial Sobolev estimate and a logarithmic energy
inequality. That mechanism is specific to the cubic problem. Indeed, a
direct energy estimate for a high--low frequency decomposition $u=v+w$ contains a
term of the form
\begin{equation}
\left\lVert w(t)\right\rVert_{L^\infty}E(v(t))^{\frac{p+1}{4}},
\label{eq:badpower}
\end{equation}
whose power of the energy is strictly larger than one when $p>3$.

\begin{remark}
	The explicit conjugation \eqref{eq:conjugation} is special to radial data
	outside a ball. It supplies both the cancellation in
	\eqref{eq:sinekernel} and the favorable boundary term in
	\eqref{eq:Morawetzidentity}. For a nonradial solution or a general
	obstacle, neither ingredient is presently available in this form.
\end{remark}

\begin{remark}
	The argument proves global well-posedness and a uniform energy bound for
	the low-frequency remainder. It does not by itself yield a global Strichartz estimate of the solution. 
\end{remark}

Several aspects distinguish the exterior problem from the usual
whole space. A low frequency truncation of data produces a
finite-energy component, but on a fixed exterior domain its
energy may grow with the cut-off scale. The resulting energy identity also
contains the superlinear factor displayed in \eqref{eq:badpower}. Moreover, the
distorted Littlewood-Paley projections for the Dirichlet Laplacian are not translation-invariant convolutions in
the radial variable: their kernels contain a reflected term, see
\eqref{eq:sinekernel}. Finally, another useful spacetime control is needed by the
Morawetz estimate
\begin{equation*}
 \int_I\int_{\Omega_R}\frac{|v(t,x)|^{p+1}}{|x|}\,\mathrm{d} x\,\mathrm{d} t, 
\end{equation*}
which is used to preclude the concentration of the solution at the origin in \cite{Bourgain1999, LiSZ:NLS, LiXZ:NLS}.

The proof combines a scaling normalization with a Fourier  truncation
decomposition. Given a large dyadic number $L$, we first rescale the problem
on $\Omega_1$ to a unit-frequency problem on $\Omega_L$. The intertwining
identity
\begin{equation*}
 P_{\leq1}^{(L)}\mathcal S_L=\mathcal S_LP_{\leq L}^{(1)}
\end{equation*}
ensures that the energy of the low-frequency data is bounded independently
of $L$ (see Lemma \ref{lem:scalingcutoff}).

 Following the argument in \cite{DodsonCubic2018, Dodson2018}, we evolve small high-frequency component as a global
solution $w$ in the critical space, and write $u=v+w$. The remainder $v$ has finite energy and
satisfies a difference equation whose leading interaction is
$p|v|^{p-1}w$. Rather than estimate this term directly, we add the
correction
\begin{equation*}
 \int_{\Omega_R}|v|^{p-1}v\,w\,\mathrm{d} x
\end{equation*}
to the energy. Its derivative cancels the dangerous time-derivative
interaction. The remaining leading terms pair $F(v)$ with
$(\partial_t,\partial_r,r^{-1})w$. A weighted Littlewood--Paley
paraproduct estimate reduces these terms to
\begin{equation*}
 E(v)^{1/2}\mathcal{M}(v)^{\frac{p-1}{p+1}}B_w(t),
 \qquad
 \mathcal{M}(v):=\int_{\Omega_R}\frac{|v|^{p+1}}r\,\mathrm{d} x,
\end{equation*}
where $B_w\in L_t^{(p+1)/2}$ is smal, see Proposition \ref{prop:smallBesov}. Young's inequality absorbs the
Morawetz contribution and leaves the term
$E(v)^{(p+1)/4}B_w^{(p+1)/2}$. The coefficient is integrable and can be
made arbitrarily small by placing the Fourier cut-off sufficiently far into
the tail; a bootstrap argument then yields a global energy bound.

Both strict inequalities in $3<p<5$ are used. The lower bound gives the
superlinear exponent $(p+1)/4>1$ and the summable low-frequency power
$2-6/(p+1)>0$. The upper bound yields the positive high-frequency gain
\begin{equation*}
 \eta_p=\frac{(p-3)(5-p)}{2(p^2-1)}>0
\end{equation*}
in the weighted Besov estimate, see \eqref{def:etap} in  Lemma \ref{lem:forcedBesov},  and provides the time-integrability
required by the perturbation argument. The relevant gains
degenerate at the endpoints, and the present argument is therefore confined
to the open interval $3<p<5$.

The energy--Morawetz estimate is inspired by Dodson in
\cite{DodsonCubic2018,Dodson2018}. In the exterior setting, however,
Euclidean Fourier multipliers must be replaced by spectral multipliers of
the Dirichlet Laplacian, and normalizing the frequency changes the radius of
the obstacle. All estimates needed below are therefore formulated uniformly
over $R>0$. Further applications of Morawetz estimates to energy-critical
Dirichlet/Neumann problems can be found in
\cite{BurqLP:NLW,  BurqP:def NLW: Neum prob, DuyLaf:NLW:Neum prob, LiSZ:NLS, LiXZ:NLS, SmithSogge95:NLW JAMS} and the references therein.  

The paper is organized as follows. Section~\ref{sec:spectral} develops the
radial distorted Fourier transform, the associated distorted Littlewood--Paley
projections, and local well-posedness and stability result.
Section~\ref{sec:weighted} establishes the boundedness of weighted Littlewood-Paley projection,
a weighted radial Sobolev estimate, and
weighted Strichartz estimates. Section~\ref{sec:decomposition} constructs the 
high--low frequency decomposition. In
Section~\ref{sec:modified} we derive the exterior Morawetz identity and the energy--Morawetz estimate for the energy solution for the difference equation with low frequency component. Section~\ref{sec:mainproof} normalizes the
truncation scale and closes the global bootstrap argument.

\medskip
\noindent\textbf{Notation.}
All functions are radial unless explicitly stated otherwise. We write
$A\lesssim B$ if $A\leq CB$ for a constant independent of the obstacle
radius and the dyadic frequency; dependence on parameters is indicated by a
subscript. We write $A \approx B$ when both $A\lesssim B$ and $B\lesssim A$
hold. Mixed spacetime norms are taken over the interval and exterior domain
displayed in the argument. Dyadic numbers always belong to $2^{\mathbb Z}$.

\section{Distorted Fourier transform and local well-posedness}
\label{sec:spectral}

This section recall the Dirichlet Sobolev spaces, the  distorted
Fourier transform, and the distorted Littlewood--Paley theory; see more details in 
\cite{KVZ16:IMRN Riesz tf,LiSZ:NLS,XuYang2024}. We then review the
Strichartz estimates and stability result used below.
\subsection{Functional spaces}

Let $-\Delta_{\Omega_R}$ denote the nonnegative Friedrichs extension of
$-\Delta$ initially defined on $C_c^\infty(\Omega_R)$. For $s\in\mathbb{R}$,
the homogeneous Dirichlet space is defined by
\begin{equation*}
 \dot H_D^s(\Omega_R)
 :=\left\{f:\left(\sqrt{-\Delta_{\Omega_R}}\right)^sf\in L^2(\Omega_R)\right\},
 \qquad
 \left\lVert f\right\rVert_{\dot H_D^s}:=\left\lVert \left(\sqrt{-\Delta_{\Omega_R}}\right)^sf\right\rVert_2,
\end{equation*}
where, for $s<0$, the right-hand side is understood by completion and
duality. More explicitly, we have
\begin{equation*}
 \dot H_D^{-s}(\Omega_R)
 =\bigl(\dot H_D^s(\Omega_R)\bigr)^*,
 \qquad
 |\left\langle f,g\right\rangle|\leq\left\lVert f\right\rVert_{\dot H_D^s}\left\lVert g\right\rVert_{\dot H_D^{-s}}.
\end{equation*}
All spaces in this paper are restricted to the radial subspace unless
otherwise specified. Smooth radial functions compactly supported away
from the boundary are dense in the norms used below. The spectral resolution of
identity in \cite{LiSZ:NLS, XuYang2024} may first be established and then passed to the
completion.

The linear wave propagation is unitary in the energy space. 
\begin{align*}
 \left\lVert \cos(t\sqrt{-\Delta_{\Omega_R}})f\right\rVert_{\dot H_D^s}
 \leq\left\lVert f\right\rVert_{\dot H_D^s},
 \quad
 \left\lVert \frac{\sin(t\sqrt{-\Delta_{\Omega_R}})}{ \sqrt{-\Delta_{\Omega_R}}} g\right\rVert_{\dot H_D^s}
 \leq\left\lVert g\right\rVert_{\dot H_D^{s-1}}. 
\end{align*}
These bounds follow from the inequalities
$|\cos(t\lambda)|\leq1$ and $|\sin(t\lambda)|\leq1$.

\subsection{The distorted Fourier transform}
Following the observation about the half-wave relations between $\mathbb{R}^3$ and $\Omega_R$ in Appendix A in
\cite{XuYang2024}, we can write
$r=\left\lvert x\right\rvert=R+\rho$, where $\rho>0$, and define for a radial function
$f$ that 
\begin{equation}
 (\mathcal U_Rf)(\rho):=(R+\rho)f(R+\rho).
 \label{eq:UR}
\end{equation}
\begin{proposition}\label{prop:spectralresolution}
Up to the factor $4\pi$, $\mathcal U_R$ is unitary from
$L^2_{\mathrm{rad}}(\Omega_R)$ onto $L^2(0,\infty)$, and
\begin{equation}
 \mathcal U_R(-\Delta_{\Omega_R})\mathcal U_R^{-1}
 =-\partial_\rho^2
 \label{eq:conjugation}
\end{equation}
with the Dirichlet condition at $\rho=0$. If $c_0=(2/\pi)^{1/2}$, then
the distorted Fourier transform and its inverse are
\begin{equation}
 \mathcal F_Rf(\lambda)
 =c_0\int_0^\infty
 \sin(\lambda\rho)(\mathcal U_Rf)(\rho)\,\mathrm{d}\rho,
 \label{eq:DFT}
\end{equation}
and
\begin{equation}
 f(R+\rho)
 =\frac{c_0}{R+\rho}\int_0^\infty
 \sin(\lambda\rho)\mathcal F_Rf(\lambda)\,\mathrm{d}\lambda.
 \label{eq:inverseDFT}
\end{equation}
Moreover, for every $s\in\mathbb{R}$, we have
\begin{equation}
 \left\lVert f\right\rVert_{\dot H^s_D(\Omega_R)}^2
 =4\pi
 \int_0^\infty\lambda^{2s}
 \left\lvert \mathcal F_Rf(\lambda)\right\rvert^2\,\mathrm{d}\lambda.
 \label{eq:spectralnorm}
\end{equation}
\end{proposition}

\begin{proof} 
The result is corresponding to the distorted (inverse) Fourier transform in Section~2 of \cite{XuYang2024}. 
For a smooth radial Dirichlet function $f=h(\rho)/(R+\rho)$, direct
differentiation gives
\begin{equation*}
 -\Delta f(R+\rho)=\frac{-h''(\rho)}{R+\rho},
 \qquad h(0)=0.
\end{equation*}
Also
$\left\lVert f\right\rVert_{L^2(\Omega_R)}^2=4\pi\left\lVert h\right\rVert_{L^2(0,\infty)}^2$.
Thus the Friedrichs realization of the radial Dirichlet Laplacian is
unitarily equivalent to the half-line Dirichlet Laplacian.  The normalized
sine-Fourier transform diagonalizes the latter and is an involutive unitary map.
 \eqref{eq:DFT}--\eqref{eq:spectralnorm} follow first for smooth
compactly supported functions and then by completion. This completes the proof.
\end{proof}

Choose $\phi\in C_c^\infty([0,2))$, with $\phi=1$ on $[0,1]$, and set
$\psi(\lambda)=\phi(\lambda)-\phi(2\lambda)$. For dyadic $N$, we define
\begin{equation*}
 \begin{aligned}
 P_{\leq N}^{(R)}&:=\phi(\sqrt{-\Delta_{\Omega_R}}/N),
 &\qquad P_N^{(R)}&:=\psi(\sqrt{-\Delta_{\Omega_R}}/N),\\
 P_{>N}^{(R)}&:=I-P_{\leq N}^{(R)},
 &P_{\geq N}^{(R)}&:=I-P_{\leq N/2}^{(R)}.
 \end{aligned}
\end{equation*}
We use $\widetilde P_N^{(R)}$ for a smooth annular multiplier that equals
one on the frequency support of $P_N^{(R)}$ and is supported where
$\lambda\approx N$.
All implicit constants below are independent of $R$ and the dyadic
frequency unless otherwise stated.

\begin{proposition}\label{prop:dyadiccalculus}
For $s\in\mathbb{R}$ and radial $f\in\dot H_D^s(\Omega_R)$, we have
\begin{equation}
 \left\lVert f\right\rVert_{\dot H_D^s}^2
 \approx_\phi\sum_{N\in2^\mathbb Z}
 N^{2s}\left\lVert P_N^{(R)}f\right\rVert_2^2.
 \label{eq:squarefunctionL2}
\end{equation}
If the frequency support of $f_N$ is contained in
$[N/4,4N]$, then for any $s_0,s_1\in\mathbb{R}$, we have
\begin{equation}
 \left\lVert f_N\right\rVert_{\dot H_D^{s_1}}
\approx N^{s_1-s_0}\left\lVert f_N\right\rVert_{\dot H_D^{s_0}}.
 \label{eq:spectralBernstein}
\end{equation}
Moreover, spectral multipliers commute with the wave propagation and with each
other:
\begin{equation}
 P_N^{(R)}m(\sqrt{-\Delta_{\Omega_R}})=m(\sqrt{-\Delta_{\Omega_R}})P_N^{(R)},\qquad
 P_N^{(R)}\mathcal W_R(t)=\mathcal W_R(t)P_N^{(R)}.
 \label{eq:multipliercommutation}
\end{equation}
For $1<r<\infty$, the Littlewood--Paley square function satisfies
\begin{equation}
 \left\lVert f\right\rVert_{L^r(\Omega_R)}
\approx_r
 \left\lVert \left(\sum_{N\in2^{\mathbb Z}}
 |P_N^{(R)}f|^2\right)^{1/2}\right\rVert_{L^r(\Omega_R)},
 \label{eq:squarefunctionLp}
\end{equation}
with constants uniform in $R$.
\end{proposition}

\begin{proof}
See \cite{KVZ16:IMRN Riesz tf,LiSZ:NLS,XuYang2024} for the general
multiplier framework. Here we recall the spectral argument. By
\eqref{eq:spectralnorm}, we have
\begin{align*}
 \sum_NN^{2s}\left\lVert P_N^{(R)}f\right\rVert_2^2
 &=4\pi\int_0^\infty
 \left(\sum_NN^{2s}|\psi(\lambda/N)|^2\right)
 |\mathcal F_Rf(\lambda)|^2\,\mathrm{d}\lambda.
\end{align*}
Only finitely many summands occur at each $\lambda>0$, and on their
support $N\approx\lambda$. The expression in $\sum_NN^{2s}|\psi(\lambda/N)|^2$ is therefore
comparable to $\lambda^{2s}$, this proves
\eqref{eq:squarefunctionL2}. The same support observation shows
\eqref{eq:spectralBernstein}. Every operator in
\eqref{eq:multipliercommutation} is multiplication by a scalar function
of $\lambda$ after applying $\mathcal F_R$, so the commutation identities
follow. 

At last, the $L^r$ statement in \eqref{eq:squarefunctionLp} is obtained from the Mikhlin
multiplier theorem for the Dirichlet Laplacian and the usual randomization
argument by the standard Khintchine inequality.  
\end{proof}

We shall repeatedly use the following two facts. For $s<1$, we have
\begin{align}
 \left\lVert P_{\leq N}^{(R)}f\right\rVert_{\dot H_D^1}
 &\lesssim N^{1-s}\left\lVert f\right\rVert_{\dot H_D^s}, \quad 
 \left\lVert P_{>N}^{(R)}f\right\rVert_{\dot H_D^s}
 \longrightarrow 0\quad \text{as}\; N\to\infty.\label{eq:spectralTail}
\end{align}
To make the second assertion explicit, the spectral identity gives
\begin{equation*}
 \left\|P_{>N}^{(R)}f\right\|_{\dot H_D^s}^2
 =4\pi\int_0^\infty
 \lambda^{2s}\left|1-\phi(\lambda/N)\right|^2
 \left|\mathcal F_Rf(\lambda)\right|^2\,\mathrm d\lambda.
\end{equation*}
The multiplier tends pointwise to zero as $N\to\infty$ and is bounded by
one. The integrable majorant is the integrand in the spectral definition of
$\left\|f\right\|_{\dot H_D^s}^2$. The dominated convergence theorem
therefore proves \eqref{eq:spectralTail}.

\begin{lemma}\label{lem:basictools}
Let $0\le\delta<1$, $m_\delta=6/(3+2\delta)$ and $0\le s\le1$. For radial Dirichlet
functions on $\Omega_R$, the following Sobolev estimates and Hardy's inequality hold
\begin{align}
 \left\lVert f/r\right\rVert_{L^2(\Omega_R)}&\le2\left\lVert \nabla f\right\rVert_{L^2(\Omega_R)},
 \label{eq:Hardy}\\
 \left\lVert f\right\rVert_{L^6(\Omega_R)}&\lesssim\left\lVert \nabla f\right\rVert_{L^2(\Omega_R)},
 \label{eq:Sobolev6}\\
 \left\lVert f\right\rVert_{L^{6/(3-2s)}(\Omega_R)}
 &\lesssim_s\left\lVert f\right\rVert_{\dot H_D^s(\Omega_R)},
 \label{eq:Sobolevgeneral}\\
 \left\lVert G\right\rVert_{\dot H_D^{-\delta}(\Omega_R)}
 &\lesssim_\delta\left\lVert G\right\rVert_{L^{m_\delta}(\Omega_R)},
 \label{eq:negativeSobolev}
\end{align}
where the constants are independent of $R$. In addition, if $s_0<s<s_1$ and
$\theta=(s-s_0)/(s_1-s_0)$, then the interpolation inequality holds
\begin{equation}
 \left\lVert f\right\rVert_{\dot H_D^s}
 \le \left\lVert f\right\rVert_{\dot H_D^{s_0}}^{1-\theta}
      \left\lVert f\right\rVert_{\dot H_D^{s_1}}^\theta.
 \label{eq:spectralinterpolation}
\end{equation}
\end{lemma}

\begin{proof}
Let $\widetilde f$ be the extension of $f$ by zero to
$\{|x|\leq R\}$. Since $f$ has zero trace on $\partial\Omega_R$, the
extension theorem for $H_0^1$ gives
\begin{equation*}
 \widetilde f\in\dot H^1(\mathbb{R}^3),\qquad
 \left\|\nabla\widetilde f\right\|_{L^2(\mathbb{R}^3)}
 =\left\|\nabla f\right\|_{L^2(\Omega_R)}.
\end{equation*}
The Euclidean inequalities applied to $\widetilde f$ give
\begin{align*}
 \int_{\Omega_R}\frac{|f(x)|^2}{|x|^2}\,\mathrm dx
 \leq4\int_{\Omega_R}|\nabla f(x)|^2\,\mathrm dx,
\quad
 \left(\int_{\Omega_R}|f(x)|^6\,\mathrm dx\right)^{1/3}
 \leq C\int_{\Omega_R}|\nabla f(x)|^2\,\mathrm dx.
\end{align*}
This proves \eqref{eq:Hardy} and \eqref{eq:Sobolev6}, with constants
independent of $R$.

Interpolating the identity embedding $L^2\hookrightarrow L^2$ with
\eqref{eq:Sobolev6} yields
\begin{equation*}
 \frac1{r_s}=\frac{1-s}{2}+\frac{s}{6}
 =\frac{3-2s}{6},
 \qquad
 \left\|f\right\|_{L^{r_s}}
 \lesssim_s\left\|f\right\|_{\dot H_D^s},
\end{equation*}
which gives \eqref{eq:Sobolevgeneral}. If
$G\in L^{m_\delta}(\Omega_R)$ and
$\varphi\in C_c^\infty(\Omega_R)$, then H\"older's inequality and
\eqref{eq:Sobolevgeneral} imply
\begin{align*}
 |\langle G,\varphi\rangle|
 &\leq\left\|G\right\|_{L^{m_\delta}}
       \left\|\varphi\right\|_{L^{m_\delta'} }
 \lesssim_\delta
 \left\|G\right\|_{L^{m_\delta}}
 \left\|\varphi\right\|_{\dot H_D^\delta},
\end{align*}
because $m_\delta'=6/(3-2\delta)$. By taking the supremum over
$\left\|\varphi\right\|_{\dot H_D^\delta}=1$, we get
\eqref{eq:negativeSobolev}. Finally, if
$\theta=(s-s_0)/(s_1-s_0)$, then we have
\begin{align*}
 \int_0^\infty\lambda^{2s}|\mathcal F_Rf(\lambda)|^2\,\mathrm d\lambda
 &\leq
 \left(\int_0^\infty\lambda^{2s_0}|\mathcal F_Rf(\lambda)|^2
 \,\mathrm d\lambda\right)^{1-\theta}
 \left(\int_0^\infty\lambda^{2s_1}|\mathcal F_Rf(\lambda)|^2
 \,\mathrm d\lambda\right)^{\theta}.
\end{align*}
By taking square roots, we have \eqref{eq:spectralinterpolation}, and complete the proof.
\end{proof}

\begin{lemma}
\label{lem:energyinterpolation}
Let $v\in\dot H_D^1(\Omega_R)\cap L^{p+1}(\Omega_R)$ and
$p+1\leq m\leq6$. Define $\theta_m\in[0,1]$ by
\begin{equation*}
 \frac1m=\frac{\theta_m}{p+1}+\frac{1-\theta_m}{6}.
\end{equation*}
and put
\begin{equation*}
 \mathcal H_p(v):=\frac12\left\lVert \nabla v\right\rVert_2^2
 +\frac1{p+1}\left\lVert v\right\rVert_{p+1}^{p+1}.
\end{equation*}
Then we have
\begin{equation}
 \left\lVert v\right\rVert_{L^m}
 \lesssim_p
 \mathcal H_p(v)^{\frac{\theta_m}{p+1}
 +\frac{1-\theta_m}{2}}.
 \label{eq:energyinterpolation}
\end{equation}
In addition,
\begin{equation}
 \left\lVert v/r\right\rVert_2\lesssim \mathcal H_p(v)^{1/2},
 \qquad
 \left\lVert v\right\rVert_6\lesssim \mathcal H_p(v)^{1/2}.
 \label{eq:energyHardySobolev}
\end{equation}
\end{lemma}

\begin{proof}
The interpolation gives
\begin{equation}
 \left\lVert v\right\rVert_m\leq\left\lVert v\right\rVert_{p+1}^{\theta_m}
 \left\lVert v\right\rVert_6^{1-\theta_m}.
 \label{eq:Lebesgueinterpolation}
\end{equation}
By the definition of $\mathcal H_p(v)$, we have
\begin{equation*}
 \left\lVert v\right\rVert_{p+1}\leq((p+1)\mathcal H_p(v))^{1/(p+1)},
 \qquad
 \left\lVert \nabla v\right\rVert_2\leq(2\mathcal H_p(v))^{1/2}.
\end{equation*}
Inserting \eqref{eq:Sobolev6} into
\eqref{eq:Lebesgueinterpolation}, we have 
\eqref{eq:energyinterpolation}; \eqref{eq:energyHardySobolev} follows
 from \eqref{eq:Hardy} and \eqref{eq:Sobolev6}. 
\end{proof}

We also use the following odd version of the one-dimensional weighted
Hardy--Sobolev inequality; compare with the classical Stein--Weiss theorem in 
\cite{SteinWeiss1958}.

\begin{lemma}\label{lem:oddHardySobolev}
Let $2\le m<\infty$, $0\le\gamma<1/m'$, and
$s=\frac12-\frac1m+\gamma$, with $0\leq s<1$. If $h$ is odd, then we have
\begin{equation*}
 \left\lVert \left\lvert \rho\right\rvert^{-\gamma}h\right\rVert_{L^m(\mathbb{R})}
 \lesssim_{m,\gamma}
 \left\lVert \left\lvert D_\rho\right\rvert^{s}h\right\rVert_{L^2(\mathbb{R})}.
\end{equation*}
\end{lemma}

\begin{proof}
The case $s=0$ is trivial, since then $m=2$ and $\gamma=0$. We therefore assume
$0<s<1$. Let  $g=|D_\rho|^sh$. Since $|D_\rho|^s$ has an even symbol, and
$g$ is odd. The Riesz-potential representation and oddness structure give, for
$\rho>0$, that
\begin{equation*}
 h(\rho)=c_s\int_0^\infty K_s(\rho,\tau)g(\tau)\,\mathrm d\tau,
 \qquad
 K_s(\rho,\tau)
 =|\rho-\tau|^{s-1}-(\rho+\tau)^{s-1}.
\end{equation*}
The mean value theorem gives the two off-diagonal estimates
\begin{align}
 |K_s(\rho,\tau)|
 \leq C_s\tau\rho^{s-2},\quad &0<\tau\leq\rho/2,
 \label{eq:oddKernelLow}\\
 |K_s(\rho,\tau)|
 \leq C_s\rho\tau^{s-2},\quad &\tau\geq2\rho,
 \label{eq:oddKernelHigh}
\end{align}
whereas, when $\rho/2<\tau<2\rho$, we have
\begin{equation*}
 |K_s(\rho,\tau)|
 \leq |\rho-\tau|^{s-1}+C_s\rho^{s-1}.
\end{equation*}

We now sum these bounds without suppressing the scale factors. Let
$J_k=(2^k,2^{k+1})$ and set
$b_j=\left\|g\right\|_{L^2(J_j)}$.

 If $j\leq k-2$, then
\eqref{eq:oddKernelLow}, and the Cauchy--Schwarz inequality on $J_j$,  imply that
\begin{align}
 \left\|\rho^{-\gamma}
   \int_{J_j}K_s(\rho,\tau)g(\tau)\,\mathrm d\tau
   \right\|_{L^m(J_k)}
 \leq
 C_{s,m}\,2^{k(s-2-\gamma+1/m)}2^{3j/2}b_j
 =C_{s,m}\,2^{-3(k-j)/2}b_j.
 \label{eq:oddLowShellBound}
\end{align}
where 
$s=\frac12-\frac1m+\gamma$. 

If $j\geq k+2$, \eqref{eq:oddKernelHigh} gives that 
\begin{align*}
 \left\|\rho^{-\gamma}
   \int_{J_j}K_s(\rho,\tau)g(\tau)\,\mathrm d\tau
   \right\|_{L^m(J_k)}
\leq
 C_{s,m}\,2^{k(1-\gamma+1/m)}2^{j(s-3/2)}b_j
 =C_{s,m}\,2^{-\vartheta(j-k)}b_j,
\end{align*}
where
$
 \vartheta:=1-\gamma+\frac1m>0.
$

For $|j-k|\leq1$, rescale $\rho=2^kx$, $\tau=2^ky$.
Young's inequality on a bounded interval, with the endpoint
Hardy--Littlewood--Sobolev inequality when $\gamma=0$, yields
\begin{equation}
 \left\|\rho^{-\gamma}
   \int_{J_j}K_s(\rho,\tau)g(\tau)\,\mathrm d\tau
   \right\|_{L^m(J_k)}
 \leq C_{s,m,\gamma}b_j.
 \label{eq:oddDiagonalShellBound}
\end{equation}
Indeed, the dilation factor in the fractional integral is
$2^{k(s+1/m-1/2)}=2^{k\gamma}$ and is canceled by the weight
$\rho^{-\gamma}$. Equations \eqref{eq:oddLowShellBound} - 
\eqref{eq:oddDiagonalShellBound} show that the sequence of
$L^m(J_k)$ norms is bounded by the convolution of $(b_j)$ with an
$\ell^1(\mathbb{Z})$ kernel. Since $m\geq2$, we have
\begin{equation*}
 \left\|(b_j)\right\|_{\ell^m}
 \leq\left\|(b_j)\right\|_{\ell^2}
 =\left\|g\right\|_{L^2(0,\infty)}.
\end{equation*}
This completes the proof.
\end{proof}

\begin{lemma}\label{lem:radialsobolev}
Let $2\le q\leq6$. If
\begin{equation}
 \frac12-\frac1q\le s\le\frac32-\frac3q,
 \label{eq:radialrange}
\end{equation}
then the radial Sobolev embedding inequality holds
\begin{equation*}
 \left\lVert f\right\rVert_{L^q(\Omega_R)}
 \lesssim_{q,s,R}
 \left\lVert f\right\rVert_{\dot H_D^s(\Omega_R)}
\end{equation*}
for radial $f$. Moreover, if $1/2<s<3/2$, then the radial decay estimate holds
\begin{equation}
 \left\lvert f(R+\rho)\right\rvert
 \lesssim_s
 \frac{\rho^{s-1/2}}{R+\rho}
 \left\lVert f\right\rVert_{\dot H_D^s(\Omega_R)}.
 \label{eq:pointwise}
\end{equation}
\end{lemma}

\begin{proof}
Let $h=\mathcal U_Rf$ and extend it oddly.

At the lower endpoint
$s=\frac12-\frac1q$, the one-dimensional Sobolev inequality  gives
\begin{equation*}
 \left\lVert f\right\rVert_{L^q(\Omega_R)}
 \lesssim R^{2/q-1}\left\lVert h\right\rVert_{L^q(0,\infty)}
 \lesssim_{q,R}\left\lVert \left\lvert D_\rho\right\rvert^s h_{\rm odd}\right\rVert_2,
\end{equation*}
where we used 
$R+\rho\ge R$ and the fact that
\begin{equation*}
 \left\|f\right\|_{L^q(\Omega_R)}^q
 =4\pi\int_0^\infty
 |h(\rho)|^q(R+\rho)^{2-q}\,\mathrm d\rho.
\end{equation*}

At the upper endpoint
$s=\frac32-\frac3q$. Suppose first that $2\leq q<6$ and  set $\gamma=1-2/q$. Then
$0\leq\gamma<1/q'$ and $0\leq s<1$. Furthermore, we have
\begin{equation*}
 (R+\rho)^{2-q}\leq\rho^{2-q},
 \qquad
 \left\|f\right\|_{L^q(\Omega_R)}
 \leq(4\pi)^{1/q}
 \left\|\rho^{-\gamma}h\right\|_{L^q(0,\infty)}.
\end{equation*}
Lemma~\ref{lem:oddHardySobolev} gives the result. If $q=6$,
the upper endpoint is $s=1$, and the result follows from
\eqref{eq:Sobolev6}.

Interpolation between the two endpoint estimates gives
\eqref{eq:radialrange}.

To prove \eqref{eq:pointwise}, we obtain  by \eqref{eq:inverseDFT} and the
Cauchy--Schwarz inequality that
\begin{align*}
 \left\lvert f(R+\rho)\right\rvert
 \lesssim
 \frac1{R+\rho}
 \left(\int_0^\infty
       \sin^2(\lambda\rho)\lambda^{-2s}\,\mathrm{d}\lambda\right)^{1/2}
 \left\lVert f\right\rVert_{\dot H_D^s}. 
\end{align*}
After the change of variables $\mu=\lambda\rho$, the integral is
\begin{equation*}
 \int_0^\infty\sin^2(\lambda\rho)\lambda^{-2s}\,\mathrm d\lambda
 =\rho^{2s-1}\int_0^\infty\sin^2\mu\,\mu^{-2s}\,\mathrm d\mu \lesssim_s  \rho^{2s-1}, 
\end{equation*}
where we use the fact that $1/2<s<3/2$ in the last inequality.
Therefore, we have
\begin{align*}
\left\lvert f(R+\rho)\right\rvert
  \lesssim_s
\frac{\rho^{s-1/2}}{R+\rho}
\left\lVert f\right\rVert_{\dot H_D^s}.
\end{align*}This proves
\eqref{eq:pointwise}.
\end{proof}
%
%Taking $s=s_c$ and $q=p+1$ in \eqref{eq:radialrange} proves
%\eqref{eq:criticalLp}; the upper inequality is equivalent to $p\le5$.

\subsection{Strichartz estimates and local well-posedness}

Recall that
\begin{equation*}
 \mathcal W_R(t)(f,g):=\cos(t\sqrt{-\Delta_{\Omega_R}})f+\frac{\sin(t\sqrt{-\Delta_{\Omega_R}})}{\sqrt{-\Delta_{\Omega_R}}} g.
\end{equation*}

We  recall the frequency-localized Strichartz estimate for the radial
Dirichlet half-wave propagator.
\begin{lemma}
\label{lem:dyadicStrichartz}
Let $2\leq Q,S\leq\infty$ and
$Q^{-1}+2S^{-1}<1$. Then
\begin{equation*}
 \left\lVert e^{it\sqrt{-\Delta_{\Omega_R}}}P_N^{(R)}f\right\rVert_{L_t^QL_x^S(\mathbb{R}\times\Omega_R)}
 \lesssim_{Q,S}
 N^{\frac32-\frac1Q-\frac3S}\left\lVert P_N^{(R)}f\right\rVert_2.
\end{equation*}
The same estimate holds for $e^{-it\sqrt{-\Delta_{\Omega_R}}}$, uniformly in $R$ and $N$.
\end{lemma}

\begin{proof}
Please refer to Proposition~3.1 of \cite{XuYang2024}.
\end{proof}

The following estimate is the Strichartz estimate of the wave equation in the exterior domain, which is analogue to the standard results in \cite{GinV:NLW:St est,KeelTao1998}. We can also refer to \cite{Burq03: NLW CPDE, IvLebeauP:LW:dis est I, IvLasLebP:LW:dis est II, SmithSogge00:NLW CPDE} and reference therein.

\begin{lemma}\label{lem:radialStrichartz}
Let $2\le Q\le\infty$ and $2\le S<\infty$ satisfy the condition that
\begin{equation}
 \frac1Q+\frac2S<1,
 \qquad
 \frac1Q+\frac3S=\frac32-s.
 \label{eq:radialadmissible}
\end{equation}
Then
\begin{equation*}
 \left\lVert \mathcal W_R(t)(f,g)\right\rVert_{L_t^QL_x^S(\mathbb{R}\times\Omega_R)}
 \lesssim_{Q,S}
 \left\lVert f\right\rVert_{\dot H_D^s}+\left\lVert g\right\rVert_{\dot H_D^{s-1}}.
\end{equation*}
If $G\in L_t^1\dot H_D^{s-1}(I\times\Omega_R)$, then we have
\begin{align}
 &\left\lVert \int_{t_0}^t\left(\sqrt{-\Delta_{\Omega_R}}\right)^{-1}\sin((t-\tau)\sqrt{-\Delta_{\Omega_R}})G(\tau)\,\mathrm{d}\tau\right\rVert_
 {L_t^\infty\dot H_D^s\cap L_t^QL_x^S(I)}\notag\\
 &\quad+
 \left\lVert \partial_t\int_{t_0}^t\left(\sqrt{-\Delta_{\Omega_R}}\right)^{-1}\sin((t-\tau)\sqrt{-\Delta_{\Omega_R}})G(\tau)\,\mathrm{d}\tau\right\rVert_
 {L_t^\infty\dot H_D^{s-1}(I)}
 \lesssim_{Q,S}\left\lVert G\right\rVert_{L_t^1\dot H_D^{s-1}(I)}.
 \label{eq:inhomogeneousStrichartz}
\end{align}
All constants are independent of $R$.
\end{lemma}

\begin{proof}
This follows from Theorem~3.3 in  \cite{XuYang2024}; see also
\cite{DuyYang24:LW Str} by an argument based on the resolvent estimate and local energy
estimates.
\end{proof}

For an interval $I$, we define
\begin{align}
 \left\lVert u\right\rVert_{S_c(I)}
 :=\ \left\lVert u\right\rVert_{L_{t,x}^{2(p-1)}(I\times\Omega_R)}
 +\left\lVert u\right\rVert_{L_t^{p-1}L_x^{3(p-1)}(I\times\Omega_R)} +\left\lVert u\right\rVert_{L_t^2L_x^{3/(1-s_c)}(I\times\Omega_R)}.
 \label{eq:Sc}
\end{align}
Each exponent pair in \eqref{eq:Sc} satisfies \eqref{eq:radialadmissible} with
$s=s_c$. We revisit the nonlinear estimate needed for the
fixed-point argument.

\begin{lemma}\label{lem:criticalnonlinear}
For any interval $I$ and radial $u,v$, we have
\begin{align}
 \left\lVert F(u)-F(v)\right\rVert_{L_t^1\dot H_D^{s_c-1}(I)} \lesssim_p
 \left(\left\lVert u\right\rVert_{L_t^{p-1}L_x^{3(p-1)}(I)}^{p-1}
 +\left\lVert v\right\rVert_{L_t^{p-1}L_x^{3(p-1)}(I)}^{p-1}\right)
 \left\lVert u-v\right\rVert_{L_t^\infty\dot H_D^{s_c}(I)}.
 \label{eq:criticaldifference}
\end{align}
In particular,
\begin{equation*}
 \left\lVert F(u)\right\rVert_{L_t^1\dot H_D^{s_c-1}(I)}
 \lesssim_p
 \left\lVert u\right\rVert_{S_c(I)}^{p-1}
 \left\lVert u\right\rVert_{L_t^\infty\dot H_D^{s_c}(I)}.
 \label{eq:chain}
\end{equation*}
\end{lemma}

\begin{proof}
It suffices to prove \eqref{eq:criticaldifference}.
Let $\delta=1-s_c$ and
$m=6/(3+2\delta)=3(p-1)/(p+1)$. By
\eqref{eq:negativeSobolev} and the difference estimate, we have
\begin{align}
 \left\lVert F(u)-F(v)\right\rVert_{\dot H_D^{s_c-1}}
 &\lesssim
 \left\lVert (\left\lvert u\right\rvert^{p-1}+\left\lvert v\right\rvert^{p-1})\left\lvert u-v\right\rvert\right\rVert_{L^m}\notag\\
 &\lesssim
 \left(\left\lVert u\right\rVert_{L^{3(p-1)}}^{p-1}
 +\left\lVert v\right\rVert_{L^{3(p-1)}}^{p-1}\right)
 \left\lVert u-v\right\rVert_{L^{3(p-1)/2}}.
 \label{eq:criticaldifferencepointwise}
\end{align}
Here the first inequality follows from the negative Sobolev estimate in
Lemma~\ref{lem:basictools}, and the second follows from H\"older's inequality,
since
\begin{equation*}
 \frac1m
 =\frac{p-1}{3(p-1)}+\frac{2}{3(p-1)}.
\end{equation*}
Thus \eqref{eq:Sobolevgeneral} at $s=s_c$ gives
\begin{equation*}
 \left\|u-v\right\|_{L^{3(p-1)/2}}
 \lesssim_p\left\|u-v\right\|_{\dot H_D^{s_c}}.
\end{equation*}
Integrating \eqref{eq:criticaldifferencepointwise} in time gives
\begin{align*}
 \int_I\left\|u(t)\right\|_{L^{3(p-1)}}^{p-1}
 \left\|u(t)-v(t)\right\|_{\dot H_D^{s_c}}\,\mathrm dt
\leq
 \left\|u\right\|_{L_t^{p-1}L_x^{3(p-1)}(I)}^{p-1}
 \left\|u-v\right\|_{L_t^\infty\dot H_D^{s_c}(I)}.
\end{align*}
This completes the proof.
\end{proof}

\begin{proposition}\label{prop:local}
Let $3<p<5$. Equation~\eqref{eq:nlw} is locally well posed in
$\dot H_D^{s_c}\times\dot H_D^{s_c-1}$. Moreover, sufficiently small
solutions are global and scatter in both time directions. More precisely, there
exists $\varepsilon_p>0$
such that, whenever
\begin{equation*}
 \left\lVert (u_0,u_1)\right\rVert_{\dot H_D^{s_c}\times\dot H_D^{s_c-1}}
 \le\varepsilon_p,
\end{equation*}
the solution of \eqref{eq:nlw} is global and
\begin{equation}
 \left\lVert u\right\rVert_{S_c(\mathbb{R})}
 +\left\lVert (u,\partial_tu)\right\rVert_{L_t^\infty(\dot H_D^{s_c}\times
 \dot H_D^{s_c-1})}
 \lesssim_p
 \left\lVert (u_0,u_1)\right\rVert_{\dot H_D^{s_c}\times\dot H_D^{s_c-1}}.
 \label{eq:smalldata}
\end{equation}
The same statement holds on every $\Omega_R$, with a constant independent
of $R$.
\end{proposition}

\begin{proof} We can refer to 
	\cite{XuYangYou2025}. 
We include the proof for completeness; For $t_0\in I$, we define the  map
\begin{equation*}
 \Phi(u)(t)=\mathcal W_R(t-t_0)(u_0,u_1)
 -\int_{t_0}^t \frac{\sin((t-s)\sqrt{-\Delta_{\Omega_R}})}{\sqrt{-\Delta_{\Omega_R}}} F(u(s))\,\mathrm{d} s.
\end{equation*}

Introduce
\begin{align*}
 \left\lVert u\right\rVert_{X(I)}
 :=&\ \left\lVert u\right\rVert_{S_c(I)}
 +\left\lVert u\right\rVert_{L_t^\infty\dot H_D^{s_c}(I)}
 +\left\lVert \partial_tu\right\rVert_{L_t^\infty\dot H_D^{s_c-1}(I)},\\
 A_I(u)
 :=&\ \left\lVert u\right\rVert_{L_t^{p-1}L_x^{3(p-1)}(I)}.
\end{align*}
Let $$K_0=\left\lVert (u_0,u_1)\right\rVert_{\dot H_D^{s_c}\times
\dot H_D^{s_c-1}}.$$
By Lemma~\ref{lem:radialStrichartz}, we have
\begin{equation}
 \left\lVert \Phi(u)\right\rVert_{X(I)}
 \leq C_pK_0+C_pA_I(u)^{p-1}
 \left\lVert u\right\rVert_{L_t^\infty\dot H_D^{s_c}(I)}.
 \label{eq:Phibound}
\end{equation}
The difference estimate gives
\begin{align}
 \left\lVert \Phi(u)-\Phi(v)\right\rVert_{X(I)}
 \leq C_p\bigl(A_I(u)^{p-1}+A_I(v)^{p-1}\bigr)
 \left\lVert u-v\right\rVert_{X(I)}.
 \label{eq:Phidifference}
\end{align}
Choose $\eta=\eta(p,K_0)>0$ so that
\begin{equation}
 2^pC_p\eta^{p-1}\leq\frac14,\qquad
 C_p(2\eta)^{p-1}(2C_pK_0+2\eta)\leq\eta.
 \label{eq:localetachoice}
\end{equation}
This is possible because $p>2$. Global integrability in time of the free wave propagation
permits a short interval $I\ni t_0$ for which
\begin{equation*}
 A_I\bigl(\mathcal W_R(t-t_0)(u_0,u_1)\bigr)\leq\eta.
\end{equation*}
On the closed subset
\begin{equation*}
 \mathfrak B_I:=\left\{u:\left\lVert u\right\rVert_{X(I)}\leq2C_pK_0+2\eta,\quad
 A_I(u)\leq2\eta\right\},
\end{equation*}
the $A_I$ component of \eqref{eq:Phibound} and
\eqref{eq:Phidifference} show that $\Phi$ is invariant and has Lipschitz
constant at most $1/2$. Banach's fixed point theorem implies a unique solution. 

If $I_{\max}=(T_-,T_+)$ is the maximal lifespan, the preceding argument
can be restarted at any $t_0\in I_{\max}$. In particular, the following
blow-up criterion holds:
\begin{equation*}
 T_+<\infty
 \quad\Longrightarrow\quad
 \left\lVert u\right\rVert_{L_t^{p-1}L_x^{3(p-1)}([0,T_+)\times\Omega_R)}=\infty,
\end{equation*}
and the analogous assertion holds backward in time.  

If $K_0\leq\varepsilon_p$, choose $\varepsilon_p$ so that the global
Strichartz norm $A_{\mathbb{R}}$ is at most $C_pK_0$ and set
$\eta=2C_pK_0$. After decreasing $\varepsilon_p$, both conditions in
\eqref{eq:localetachoice} hold. The preceding contraction then works on
$\mathbb{R}$, and absorption in \eqref{eq:Phibound} proves
\eqref{eq:smalldata}.  This completes the proof.
\end{proof}

\begin{corollary}\label{cor:stability}
Let the approximate solution $\widetilde u$ satisfy
\begin{equation*}
\left\{
\begin{aligned}
&\partial_t^2\widetilde u-\Delta \widetilde u+F(\widetilde u)=e,
&&(t,x)\in I\times\Omega_R,\\
&\widetilde u(t,x)=0,&&x\in\partial\Omega_R,\\
&(\widetilde u,\partial_t \widetilde u)|_{t=t_0}=(\widetilde u_0,\widetilde u_1),
\end{aligned}
\right.
\qquad
\end{equation*}
on an interval $I$, and suppose $A_I(\widetilde u)<\infty$. If $u$ is the solution of \eqref{eq:nlw} with data $(u_0, u_1)$
at $t_0\in I$, then
\begin{align}
 \left\lVert u-\widetilde u\right\rVert_{X(I)}
 \lesssim_p{}&
 \left\lVert (u-\widetilde u,\partial_tu-\partial_t\widetilde u)(t_0)\right\rVert
_ {\dot H_D^{s_c}\times\dot H_D^{s_c-1}} +\left\lVert e\right\rVert_{L_t^1\dot H_D^{s_c-1}(I)},\label{eq:stabilityestimate}
\end{align}
provided the right-hand side is sufficiently small.
\end{corollary}

\begin{proof}
The argument is standard; see \cite{CKSTT2008, KenigM2006:NLS, KenigM2008:NLW, KillpV2008, TaoV2005} . We provide the
details needed here. Let $z=u-\widetilde u$. Subtracting the two equations gives
\begin{equation*}
 \partial_t^2z-\Delta z
 =-\bigl(F(u)-F(\widetilde u)\bigr)-e.
\end{equation*}
The Duhamel formula is
\begin{align*}
 z(t)
 =\mathcal W_R(t-t_0)\bigl(z(t_0),\partial_tz(t_0)\bigr)
 -\int_{t_0}^t \frac{ \sin\!\left((t-s)\sqrt{-\Delta_{\Omega_R}}\right)}{\sqrt{-\Delta_{\Omega_R}}} 
 \bigl(F(u)-F(\widetilde u)+e\bigr)(s)\,\mathrm ds.
\end{align*}
By \eqref{eq:inhomogeneousStrichartz} and
\eqref{eq:criticaldifference}, we have
\begin{equation}
 \left\lVert u-\widetilde u\right\rVert_{X(I)}
 \leq C_p d_0+
 C_p\bigl(A_I(u)^{p-1}+A_I(\widetilde u)^{p-1}\bigr)
 \left\lVert u-\widetilde u\right\rVert_{X(I)},
 \label{eq:stabilityabsorption}
\end{equation}
where $d_0$ denotes the right-hand side of
\eqref{eq:stabilityestimate}. 

First suppose that $A_I(\widetilde u)\leq\eta$, where $\eta>0$ will be
chosen sufficiently small. For $t\geq t_0$, define
\begin{equation*}
 I_t=I\cap[t_0,t],
 \qquad
 \mathcal T=\left\{t:
 A_{I_t}(u)\leq2\eta,
 \ \left\|z\right\|_{X(I_t)}\leq2C_pd_0\right\}.
\end{equation*}
The set $\mathcal T$ is nonempty and closed.

On $\mathcal T$,
$A_{I_t}(u)\leq A_{I_t}(\widetilde u)+A_{I_t}(z)$ and
\eqref{eq:stabilityabsorption} imply, after choosing the constant so that
\begin{equation*}
 C_p\bigl((2\eta)^{p-1}+\eta^{p-1}\bigr)\leq\frac14,
 \qquad 2C_pd_0\leq\frac\eta2,
\end{equation*}
this implies the refined estimates
\begin{equation*}
 \left\|z\right\|_{X(I_t)}\leq\frac43C_pd_0,
 \qquad
 A_{I_t}(u)\leq\frac32\eta.
\end{equation*}
Hence $\mathcal T$ is open and equals the forward part of $I$; the
backward part is the same argument. This proves \eqref{eq:stabilityestimate}.

If $A_I(\widetilde u)<\infty$, absolute continuity of the
$L_t^{p-1}$ norm provides a partition
$I=I_1\cup\cdots\cup I_J$ such that
$A_{I_j}(\widetilde u)\leq\eta$ for every $j$. Applying the short-time
stability successively, with the terminal data on $I_j$ used as the initial
data on $I_{j+1}$, yields a finite-interval stability constant depending
only on $p$, $J$, and $A_I(\widetilde u)$. This completes the proof.
\end{proof}

\subsection{Usefule global control}
We recall an useful inequality in this part, which will be used in the global energy-Morawetz estimate of the solution of the difference equation with the low frequency part at the end of the paper. 
\begin{lemma}
	\label{lem:Bihari}
	Let $\alpha>1$, and let $y>0$ be absolutely continuous on $[0,T]$ with
	\begin{equation*}
	y'(t)\le a(t)y(t)+b(t)y(t)^\alpha,
	\qquad a,b\ge0,\quad a,b\in L^1(0,T).
	\end{equation*}
	If $A_T=\int_0^Ta(s)\,\mathrm{d} s$, then we have
	\begin{align}
	y(t)^{1-\alpha}
	\ge &\;  e^{-(\alpha-1)A_T}
	\left(y(0)^{1-\alpha}
	-(\alpha-1)e^{(\alpha-1)A_T}\int_0^tb(s)\,\mathrm{d} s\right) \notag\\
	= &\;  e^{-(\alpha-1)A_T}
	y(0)^{1-\alpha}
	-(\alpha-1)\int_0^tb(s)\,\mathrm{d} s. 
	\label{eq:Biharilem}
	\end{align}
	In particular, $y$ remains bounded whenever the right hand side is positive.
\end{lemma}

\begin{proof}
Define $z(t)=e^{-\int_0^ta(s)\,\mathrm{d} s}y(t)$,  then we have
	$$z'\le e^{(\alpha-1)\int_0^ta}b\,z^\alpha. $$ Since $1-\alpha<0$, we get
	\begin{equation*}
	\frac{\,\mathrm{d}}{\,\mathrm{d} t}z^{1-\alpha}
	\ge-(\alpha-1)e^{(\alpha-1)A_T}b(t).
	\end{equation*}
	By integration and
	$y^{1-\alpha}=e^{-(\alpha-1)\int_0^ta}z^{1-\alpha}$, we have
	\eqref{eq:Biharilem}, and complete the proof.
\end{proof}

\section{Global weighted Strichartz estimates of solution with small high-frequency data}
\label{sec:weighted}
In this section, we will show global weighted Strichartz estimates of solution of \eqref{eq:nlw} with small high-frequency data
$P_{>1}^{(R)}(u_0,u_1)$ in Proposition \ref{prop:smallBesov}.

\begin{lemma}
\label{lem:weightedLP}
For $1<q<\infty$, radial $f$, $R>0$, and dyadic $N$, we have
\begin{equation*}
 \left\lVert r^{-1/q}P_{ N}^{(R)}f\right\rVert_{L^q(\Omega_R)}+ \left\lVert r^{-1/q}P_{\le N}^{(R)}f\right\rVert_{L^q(\Omega_R)}
 +\left\lVert r^{-1/q}P_{>N}^{(R)}f\right\rVert_{L^q(\Omega_R)}
 \lesssim_q
 \left\lVert r^{-1/q}f\right\rVert_{L^q(\Omega_R)}.
\end{equation*}
\end{lemma}

\begin{proof}
Let $h=\mathcal U_Rf$. By \eqref{eq:UR}, we have
\begin{equation}
 \left\lVert r^{-1/q}f\right\rVert_{L^q(\Omega_R)}^q
 =4\pi\int_0^\infty
 (R+\rho)^{1-q}\left\lvert h(\rho)\right\rvert^q\,\mathrm{d}\rho.
 \label{eq:weightedh}
\end{equation}
The kernel of
$\mathcal U_RP_N^{(R)}\mathcal U_R^{-1}$ is
\begin{align}
 K_N(\rho,\tau) 
 = & \int^{\infty}_0 \sin(\lambda \rho) \sin(\lambda \tau) \Big(\phi_N(\lambda)-\phi_{N/2}(\lambda)\Big)  \mathrm{d}\lambda\notag\\
 = & \, N\bigl[k(N(\rho-\tau))-k(N(\rho+\tau))\bigr],
 \label{eq:sinekernel}
\end{align}
where $k$ is a Schwartz function depending only on the function $\phi$. Set
$\theta=(q-1)/q\in(0,1)$. After conjugating by the weight in
\eqref{eq:weightedh}, it suffices by Schur's lemma to prove
\begin{align}
 &\sup_{\rho>0}\int_0^\infty
 \left\lvert K_N(\rho,\tau)\right\rvert
 \left(\frac{R+\tau}{R+\rho}\right)^\theta\,\mathrm{d}\tau\lesssim_q1,
 \label{eq:schur1}\\
 &\sup_{\tau>0}\int_0^\infty
 \left\lvert K_N(\rho,\tau)\right\rvert
 \left(\frac{R+\tau}{R+\rho}\right)^\theta\,\mathrm{d}\rho\lesssim_q1.
 \label{eq:schur2}
\end{align}
By the change of variables $x=N\rho$, $y=N\tau$, and $b=NR$. In the region $y\le2x$, the
weight function obeys
\begin{equation*}
 \left(\frac{b+y}{b+x}\right)^\theta\leq3^\theta,
 \qquad 0<y\leq2x,
\end{equation*}
and the two Schwartz kernels have uniformly bounded integrals. This gives
the contribution to \eqref{eq:schur1} on that region. In the region $y>2x$, the
cancellation between the two terms in
\eqref{eq:sinekernel} and the mean value theorem give for every $M$ that
\begin{equation*}
 \left\lvert k(x-y)-k(x+y)\right\rvert
 \lesssim_M x\langle y\rangle^{-M}.
\end{equation*}
Furthermore,
\begin{equation*}
 \left(\frac{b+y}{b+x}\right)^\theta
 \le C_\theta\left(\frac yx\right)^\theta,
 \quad y>2x.
\end{equation*}
Thus it is bounded in this region by
\begin{equation*}
 Cx^{1-\theta}\int_{2x}^\infty
 y^\theta\langle y\rangle^{-M}\,\mathrm{d} y\lesssim1
\end{equation*}
when $M$ is large. This proves \eqref{eq:schur1}; the proof of
\eqref{eq:schur2} uses the complementary split $x\leq2y$ and $x>2y$.
In the latter region the same cancellation yields
$|k(x-y)-k(x+y)|\lesssim_M y\langle x\rangle^{-M}$, and
\begin{equation*}
 \left(\frac{b+y}{b+x}\right)^\theta\leq1,
 \qquad x>2y.
\end{equation*}
Hence the second Schur integral is uniform as well. The argument applies
to each smooth annular multiplier.

For $P_{\leq N}^{(R)}$, we may use 
the Schwartz kernel of $\phi(\sqrt{-\Delta_{\Omega_R}}/N)$, which has the same difference form
\eqref{eq:sinekernel}. Finally,
$P_{>N}^{(R)}=I-P_{\leq N}^{(R)}$, so  two weighted boundednesses of the Littlewood-Paley projection follow.
\end{proof}
\begin{remark}
The weighted form in Lemma~\ref{lem:weightedLP} coincides with the Morawetz
form when $q=p+1$:
\begin{equation}
 \left\lVert r^{-1/q}f\right\rVert_{L^q}^q
 =\int_{\Omega_R}\frac{\left\lvert f(x)\right\rvert^{p+1}}{\left\lvert x\right\rvert}\,\mathrm{d} x.
 \label{eq:morawetznorm}
\end{equation}
\end{remark}

\begin{lemma}\label{lem:weightedSobolev}
Assume 
\begin{equation}
q:=p+1\in(4,6),
\quad
\beta:=\frac{q-4}{2q}>0,
\quad
\sigma:=1-\frac4q, 
\label{eq:parameters}
\end{equation} we have the weighted Sobolev inequality
\begin{equation}
 \left\lVert r^\beta f\right\rVert_{L^{q/2}(\Omega_R)}
 \lesssim_p
 \left\lVert f\right\rVert_{\dot H_D^\sigma(\Omega_R)}
 \label{eq:weightedSobolev}
\end{equation}
for every radial $f$. Consequently,
\begin{equation}
 \left\lVert r^\beta P_{\ge N}^{(R)}f\right\rVert_{L^{q/2}(\Omega_R)}
 \lesssim_p
 N^{-4/q}\left\lVert \nabla f\right\rVert_{L^2(\Omega_R)}.
 \label{eq:weightedBernstein}
\end{equation}
\end{lemma}

\begin{proof}
Let $h=\mathcal U_Rf$ and let $h_{\rm odd}$ be its odd extension to $\mathbb{R}$.
A direct calculation gives
\begin{equation*}
 \left\lVert r^\beta f\right\rVert_{L^{q/2}(\Omega_R)}
 \approx
 \left\lVert (R+\rho)^{-\beta}h(\rho)\right\rVert_{L^{q/2}(0,\infty)}.
\end{equation*}
Since $R+\rho\ge\rho$, Lemma~\ref{lem:oddHardySobolev}, with
$m=q/2$, $\gamma=\beta$, and $s=\sigma$, yields
\begin{equation*}
 \left\lVert \rho^{-\beta}h\right\rVert_{L^{q/2}(0,\infty)}
 \lesssim_p
 \left\lVert \left\lvert D_\rho\right\rvert^{\sigma}h_{\rm odd}\right\rVert_{L^2(\mathbb{R})}.
\end{equation*}
%The exponents satisfy
%\begin{equation*}
% \beta-\frac2q=\sigma-\frac12.
%\end{equation*}
This proves \eqref{eq:weightedSobolev} by \eqref{eq:conjugation}. On the
support of $P_{\ge N}^{(R)}$, we have
\begin{equation*}
 \left\lVert P_{\ge N}^{(R)}f\right\rVert_{\dot H_D^\sigma}
 \le N^{\sigma-1}\left\lVert \nabla f\right\rVert_2
 =N^{-4/q}\left\lVert \nabla f\right\rVert_2,
\end{equation*}
which proves \eqref{eq:weightedBernstein}.
\end{proof}

\begin{lemma}\label{lem:weightedStrichartz}
Let $2<b<\infty$. For every radial $g\in L^2(\Omega_R)$, the following
frequency-localized weighted Strichartz estimate holds:
\begin{equation}
 \left\lVert r^{1/2}e^{it\sqrt{-\Delta_{\Omega_R}}}P_N^{(R)}g\right\rVert_{L_t^bL_r^\infty(\mathbb{R}\times(R,\infty))}
 \lesssim_b
 N^{1-1/b}\left\lVert P_N^{(R)}g\right\rVert_{L^2(\Omega_R)}.
 \label{eq:weightedStrichartz}
\end{equation}
The constant is independent of $R$ and $N$.
\end{lemma}

\begin{proof}
Let $\rho=r-R$. Formula \eqref{eq:inverseDFT} gives
\begin{equation*}
 e^{it\sqrt{-\Delta_{\Omega_R}}}P_N^{(R)}g(r)
 =\frac{c_0}{r}\int_0^\infty
 e^{it\lambda}\sin(\lambda\rho)
 \psi(\lambda/N)\mathcal F_Rg(\lambda)\,\mathrm{d}\lambda.
\end{equation*}

Let 
\begin{equation*}
 H_N(s):=\frac{c_0}{2i}\int_0^\infty
 e^{is\lambda}\psi(\lambda/N)\mathcal F_Rg(\lambda)\,\mathrm{d}\lambda.
\end{equation*}
then we have
\begin{equation*}
 r^{1/2}e^{it\sqrt{-\Delta_{\Omega_R}}}P_N^{(R)}g(r)
 =(R+\rho)^{-1/2}\bigl(H_N(t+\rho)-H_N(t-\rho)\bigr),
\end{equation*}
and Plancherel together with \eqref{eq:spectralnorm} gives
\begin{equation}
 \left\lVert H_N\right\rVert_{L^2(\mathbb{R})}\lesssim\left\lVert P_N^{(R)}g\right\rVert_{L^2(\Omega_R)}.
 \label{eq:profileL2}
\end{equation}

Set
\begin{equation*}
 K(\tau):=N^{-1/2}H_N(\tau/N),
 \qquad \varrho:=NR,
\end{equation*}
so that $K$ has Fourier support in a fixed compact subinterval of
$(0,\infty)$ and $\left\lVert K\right\rVert_2=\left\lVert H_N\right\rVert_2$.  The changes of variables
$\tau=Nt$ and $\sigma=N\rho$ show that the left-hand side of
\eqref{eq:weightedStrichartz} equals
\begin{equation}
 N^{1-1/b}
 \left\lVert \sup_{\sigma>0}(\varrho+\sigma)^{-1/2}
 \left\lvert K(\tau+\sigma)-K(\tau-\sigma)\right\rvert\right\rVert_{L_\tau^b}.
 \label{eq:mixednormrescaling}
\end{equation}
It therefore remains to prove, uniformly in $\varrho\ge0$, that
\begin{equation}
 \left\lVert \sup_{\sigma>0}(\varrho+\sigma)^{-1/2}
 \left\lvert K(\tau+\sigma)-K(\tau-\sigma)\right\rvert\right\rVert_{L_\tau^b}
 \lesssim_b\left\lVert K\right\rVert_2.
 \label{eq:unitmaximalestimate}
\end{equation}

Since $\varrho+\sigma\ge\sigma$, we split the supremum in
\eqref{eq:unitmaximalestimate} at $\sigma=1$. 

For $0<\sigma\le1$, the
fundamental theorem of calculus gives
\begin{equation*}
 \sigma^{-1/2}\left\lvert K(\tau+\sigma)-K(\tau-\sigma)\right\rvert
 \le 2\sigma^{1/2}\mathbf M(K')(\tau)
 \le 2\mathbf M(K')(\tau),
\end{equation*}
where $\mathbf M$ is the one-dimensional Hardy--Littlewood maximal
operator.  The $L^b$ boundedness of the operator $\mathbf M$  and the Bernstein estimate give
\begin{equation}
 \left\lVert \mathbf M(K')\right\rVert_{L^b}
 \lesssim_b\left\lVert K'\right\rVert_{L^b}
 \lesssim_b\left\lVert K\right\rVert_2.
 \label{eq:smalltranslationbound}
\end{equation}

For $\sigma\ge1$, both translated terms are controlled by the Peetre
maximal function
\begin{equation*}
 K_{1/2}^*(\tau):=
 \sup_{s\in\mathbb{R}}\frac{\left\lvert K(s)\right\rvert}{(1+\left\lvert s-\tau\right\rvert)^{1/2}}.
 \label{eq:Peetremaximal}
\end{equation*}
Choose $d$ with $2<d<b$.  The standard band-limited Peetre inequality,
followed by the maximal theorem, yields
\begin{align}
 K_{1/2}^*(\tau)
 &\lesssim_d \bigl[\mathbf M(\left\lvert K\right\rvert^d)(\tau)\bigr]^{1/d},
 \notag\\
 \left\lVert K_{1/2}^*\right\rVert_{L^b}
 &\lesssim_{b,d}
 \left\lVert \mathbf M(\left\lvert K\right\rvert^d)\right\rVert_{L^{b/d}}^{1/d}
 \lesssim_{b,d}\left\lVert K\right\rVert_{L^b}
 \lesssim_b\left\lVert K\right\rVert_2.
 \label{eq:Peetrebound}
\end{align}
 Equations \eqref{eq:smalltranslationbound} and
\eqref{eq:Peetrebound} prove \eqref{eq:unitmaximalestimate}.  Combining
this with \eqref{eq:mixednormrescaling} and \eqref{eq:profileL2} proves
\eqref{eq:weightedStrichartz}. 
\end{proof}

\begin{remark}
The restriction $b>2$ is necessary.  Indeed, choose a nonzero
fixed-annulus profile $K$ with $K(0)\neq0$.  For $\tau\to+\infty$, take
$\sigma=\tau$ in \eqref{eq:unitmaximalestimate}.  Since $K$ is rapidly
decreasing,
\begin{equation*}
 \sigma^{-1/2}\left\lvert K(\tau+\sigma)-K(\tau-\sigma)\right\rvert
 =\tau^{-1/2}\left\lvert K(2\tau)-K(0)\right\rvert
 \gtrsim_K\tau^{-1/2}.
\end{equation*}
The right-hand side does not belong to $L^2_\tau$.  Thus the strong
$L_t^2L_r^\infty$ estimate is false (the corresponding endpoint is at
most of weak type).  This causes no loss below, because the only
application uses $b=q/2=(p+1)/2>2$.
\end{remark}

For a radial function $z=z(x)$, we define the weighted Besov norm
\begin{equation}\label{def:Bz}
 B_z:=\sum_{N\in2^\mathbb Z}N^{-4/q}
 \left\lVert r^{1/2}P_N^{(R)}z\right\rVert_{L_r^\infty}.
\end{equation}

\begin{lemma}
\label{lem:spectralchainHardy}
Let $u=P_{\le N}^{(R)}v$, where $v$ is radial with finite energy, and let
$g_N$ have  frequency support in $[N/4,4N]$. Then we have
\begin{equation}
 \left\lvert \left\langle \sqrt{-\Delta_{\Omega_R}} F(u),g_N\right\rangle\right\rvert
 \lesssim_p
 N^\sigma\left\lVert \nabla v\right\rVert_2
 \mathcal{M}(v)^{(p-1)/q}
 \left\lVert r^{1/2}g_N\right\rVert_{L_r^\infty}.
 \label{eq:weightedchain}
\end{equation}
Moreover, if $f(R)=0$, then
we have
\begin{align}
 &\left\lVert r^{1/2}P_N^{(R)}\partial_rf\right\rVert_{L_r^\infty}
 +\left\lVert r^{1/2}P_N^{(R)}(f/r)\right\rVert_{L_r^\infty} \lesssim
 \sum_{M\in2^\mathbb Z}
 \omega(N/M)
 \left\lVert r^{1/2}P_M^{(R)}\sqrt{-\Delta_{\Omega_R}} f\right\rVert_{L_r^\infty},
 \label{eq:localizedHardy}
\end{align}
where $
\omega(s):=\min\{s,s^{1/2}\}$, $s>0
$, which implies that
\begin{equation}
 B_{\partial_rf}+B_{f/r}\lesssim B_{\sqrt{-\Delta_{\Omega_R}} f}.
 \label{eq:radialHardyBproved}
\end{equation}
\end{lemma}

\begin{proof}
All computations can be performed for smooth radial functions
with compact support, the general case follows by frequency
truncation and density. Throughout the proof, we take\begin{equation*}
q:=p+1\in(4,6),
\quad
\beta:=\frac{q-4}{2q}>0,
\quad
\sigma:=1-\frac4q, 
\end{equation*}
and
\begin{equation}
 a:=\frac q2,\qquad a'=\frac q{q-2}=\frac q{p-1},
 \qquad \frac1a+\frac1{a'}=1.
 \label{eq:chainHolderexponents}
\end{equation}

\smallskip
\noindent\emph{Step 1: two  one-dimensional  weighted estimates.}
We first establish the weighted bounds that will be used in the
nonlinear estimate. Put
\begin{equation*}
\gamma:=-\beta a=-\frac{q-4}{4}\in\left(-\frac12,0\right),
\quad
w_R(x):=(R+\lvert x\rvert)^\gamma,\quad R\geq0.
\end{equation*}

We first verify that the $A_a$-weight constant of $w_R$ is bounded
independently of $R$.  Let $I\subset\mathbb R$ be an interval, let
$\ell:=\lvert I\rvert$, and let $d:=\operatorname{dist}(I,0)$.

 If
$d\geq2\ell$, then $R+\lvert x\rvert\approx R+d$ on $I$, and therefore
\begin{equation*}
\left(\frac1\ell\int_Iw_R(x)\,\mathrm{d}x\right)
\left(\frac1\ell\int_Iw_R(x)^{-1/(a-1)}\,\mathrm{d}x\right)^{a-1}
\leq C_{a,\gamma}.
\end{equation*}

If $d<2\ell$, then $I\subset[-3\ell,3\ell]$. Since
$-1<\gamma<0$ and $-\gamma/(a-1)>0$, direct integration gives
\begin{align*}
\frac1\ell\int_I(R+\lvert x\rvert)^\gamma\,\mathrm{d}x
\leq C_\gamma(R+\ell)^\gamma,\quad 
\frac1\ell\int_I(R+\lvert x\rvert)^{-\gamma/(a-1)}\,\mathrm{d}x
\leq C_{a,\gamma}(R+\ell)^{-\gamma/(a-1)}.
\end{align*}
Multiplying these two inequalities, we have
\begin{equation*}
[w_R]_{A_a}
:=\sup_I
\left(\frac1{\lvert I\rvert}\int_Iw_R\right)
\left(\frac1{\lvert I\rvert}\int_Iw_R^{-1/(a-1)}\right)^{a-1}
\leq C_{a,\gamma}.
\end{equation*}
Therefore the weighted Hilbert-transform theorem  yields
\begin{equation*}
\left\lVert (R+\lvert\cdot\rvert)^{-\beta}
\mathcal H g\right\rVert_{L^a(\mathbb R)}
\leq C_p
\left\lVert (R+\lvert\cdot\rvert)^{-\beta}
g\right\rVert_{L^a(\mathbb R)},
\end{equation*}
with a constant independent of $R$.

We also need a one-sided Hardy inequality. Let
$g\in C_c^\infty([0,\infty))$, $g\geq0$, and
$G(\rho):=\int_0^\rho g(s)\,\mathrm{d}s$. Put
$\lambda:=a(\beta+1)>1$. Integration by parts and H\"older's inequality
give
\begin{align*}
I
\; :=& \int_0^\infty(R+\rho)^{-\lambda}G(\rho)^a\,\mathrm{d}\rho\\
=&\frac{a}{\lambda-1}
\int_0^\infty
(R+\rho)^{1-\lambda}G(\rho)^{a-1}g(\rho)\,\mathrm{d}\rho\\
\leq& \frac{a}{\lambda-1}
I^{(a-1)/a}
\left(\int_0^\infty
(R+\rho)^{-\beta a}g(\rho)^a\,\mathrm{d}\rho\right)^{1/a}.
\end{align*}
The boundary terms vanish at infinity because
$1-\lambda<0$. At $\rho=0$ they vanish because $G(0)=0$; when $R=0$,
one additionally uses
$1-\lambda+a=1-a\beta>0$. Dividing by $I^{(a-1)/a}$, with the trvial case
$I=0$ being immediate, we have 
\begin{equation*}
\left\lVert (R+\rho)^{-\beta-1}G\right\rVert_{L^a(0,\infty)}
\leq\frac{a}{a(\beta+1)-1}
\left\lVert (R+\rho)^{-\beta}g\right\rVert_{L^a(0,\infty)}.
\end{equation*}

\smallskip
\noindent\emph{Step 2: the weighted nonlinear estimate.} For $\rho>0$ and $r=R+\rho$, define 
\begin{equation*}
H(\rho):=rU(r),\quad
A(\rho):=rF(U(r)),\quad
\Gamma(\rho):=rG(r),\quad
B(\rho):=\lvert D_\rho\rvert H(\rho).
\end{equation*}
Extend all four functions oddly to $\mathbb R$. By the conjugation
formula, we have
\begin{equation*}
B=\mathcal U_R\sqrt{-\Delta_{\Omega_R}}U.
\end{equation*}
Put $Z:=\left\lVert r^{1/2}G\right\rVert_{L^\infty}$. By
the radial measure $\mathrm{d}x=4\pi r^2\,\mathrm{d}r$ and from
$B=r\sqrt{-\Delta_{\Omega_R}}U$, we  obtain $
\lvert\Gamma(\rho)\rvert
 \leq Zr^{1/2},$ and 
\begin{equation}
\begin{aligned}
\left\lVert r^{(p-1)/q}\lvert U\rvert^{p-1}
\right\rVert_{L^{a'}_\rho}^{a'}
&=\int_0^\infty r\lvert U\rvert^q\,\mathrm{d}\rho
=\frac1{4\pi}\left\lVert r^{-1/q}U
\right\rVert_{L^q(\Omega_R)}^q,\\
\left\lVert r^{-\beta}B\right\rVert_{L^a_\rho}^a
&=\int_0^\infty r^{a(1-\beta)}
\left\lvert\sqrt{-\Delta_{\Omega_R}}U\right\rvert^a\,\mathrm{d}\rho
=\frac1{4\pi}\left\lVert r^\beta
\sqrt{-\Delta_{\Omega_R}}U\right\rVert_{L^a(\Omega_R)}^a.
\end{aligned}
\label{eq:chainhalflineidentities}
\end{equation}
where we use the fact that $(p-1)a'=q$, $(p-1)a'/q=1$, and
$a(1-\beta)=\beta a+2$. 

We prove that, whenever $U$ is radial and Dirichlet and $G$ has
frequency  support in $[\Lambda/4,4\Lambda]$,
\begin{align}
\left\lvert
\left\langle\sqrt{-\Delta_{\Omega_R}}F(U),G\right\rangle
\right\rvert = &  \left\lvert  4\pi\int_0^\infty
\left\lvert D_\rho\right\rvert A(\rho)\Gamma(\rho)\,\mathrm{d}\rho \right\rvert \notag \\
\leq &\; C_p
\left\lVert r^{-1/q}U\right\rVert_{L^q(\Omega_R)}^{p-1}
\left\lVert r^\beta\sqrt{-\Delta_{\Omega_R}}U
\right\rVert_{L^a(\Omega_R)}
\left\lVert r^{1/2}G\right\rVert_{L^\infty(\Omega_R)}.
\label{eq:weightedchainkernelclaim}
\end{align}
The constant will be independent of $R$ and $\Lambda$.

The direct and reflected parts of the half-line kernel are
\begin{equation}
\lvert D_\rho\rvert A(\rho)
=c\,\operatorname{p.v.}\!\int_0^\infty
\left[
\frac{A(\rho)-A(\tau)}{(\rho-\tau)^2}
+\frac{A(\rho)+A(\tau)}{(\rho+\tau)^2}
\right]\,\mathrm{d}\tau.
\label{eq:oddhalflineDkernel}
\end{equation}

Let $\widetilde\Pi_\Lambda$ be a smooth annular multiplier equal to one
on $[\Lambda/4,4\Lambda]$, and put
$$T_\Lambda:=\mathcal H\widetilde\Pi_\Lambda.$$ Its convolution
kernel has the form
\begin{equation*}
L_\Lambda(z)=\Lambda L(\Lambda z),\qquad L\in\mathcal S(\mathbb R).
\end{equation*}

We first claim that
\begin{equation}
\sup_{\rho>0}(R+\rho)^{-1/2}
\left\lvert T_\Lambda\Gamma(\rho)-T_\Lambda\Gamma(0)\right\rvert
\leq C
\sup_{\tau>0}(R+\tau)^{-1/2}\lvert\Gamma(\tau)\rvert.
\label{eq:chainannularboundary}
\end{equation}

Denote the supremum on the right-hand side of
\eqref{eq:chainannularboundary} by $Z$. Since $\Gamma$ is odd, we have
\begin{align*}
T_\Lambda\Gamma(\rho)-T_\Lambda\Gamma(0)
&=\int_0^\infty
\bigl[L_\Lambda(\rho-\tau)-L_\Lambda(-\tau)\bigr]
\Gamma(\tau)\,\mathrm{d}\tau -\int_0^\infty
\bigl[L_\Lambda(\rho+\tau)-L_\Lambda(\tau)\bigr]
\Gamma(\tau)\,\mathrm{d}\tau.
\end{align*}
The first term in the right hand is the regularized direct term and the second term is the
regularized reflected term. Set
\begin{align*}
x&:=\Lambda\rho,\qquad y:=\Lambda\tau,\qquad b:=\Lambda R,\\
E_<&:=(0,x/2],\qquad
E_\sim:=(x/2,2x),\qquad E_>:=[2x,\infty).
\end{align*}
For $j\in\{<,\sim,>\}$ define
\begin{align*}
J_j^-(x,b)
&:=\int_{E_j}
\lvert L(x-y)-L(-y)\rvert(b+y)^{1/2}\,\mathrm{d}y,\\
J_j^+(x,b)
&:=\int_{E_j}
\lvert L(x+y)-L(y)\rvert(b+y)^{1/2}\,\mathrm{d}y.
\end{align*}
Then
\begin{equation*}
\lvert T_\Lambda\Gamma(\rho)-T_\Lambda\Gamma(0)\rvert
\leq Z\Lambda^{-1/2}
\sum_{j\in\{<,\sim,>\}}\bigl(J_j^-(x,b)+J_j^+(x,b)\bigr).
\end{equation*}
We estimate all six terms.

If $0<x\leq1$, the mean-value theorem gives, for either sign,
\begin{align*}
J_j^\pm(x,b)
&\leq x\int_0^1\int_{E_j}
\lvert L'(\theta x\pm y)\rvert(b+y)^{1/2}
\,\mathrm{d}y\,\mathrm{d}\theta\\
&\leq Cx(b+x+1)^{1/2}
\leq C(b+x)^{1/2}.
\end{align*}
For the direct term the argument of $L'$ is $\theta x-y$; for the
reflected term it is $\theta x+y$. The last inequality follows from
\begin{equation*}
x^2(b+x+1)\leq b+2x\leq2(b+x).
\end{equation*}

It remains to consider $x\geq1$. On $E_<$, one has
$\lvert x-y\rvert\geq x/2$ and $x+y\geq x$. Hence, for any $m>3$, we have
\begin{align*}
J_<^-(x,b)+J_<^+(x,b)
&\leq
C_m\int_0^{x/2}\langle x-y\rangle^{-m}(b+y)^{1/2}\,\mathrm{d}y\\
&\quad+
C_m\int_0^{x/2}\langle x+y\rangle^{-m}(b+y)^{1/2}\,\mathrm{d}y +
2\int_0^{x/2}\lvert L(y)\rvert(b+y)^{1/2}\,\mathrm{d}y\\
&\leq C_mx^{1-m}(b+x)^{1/2}+C(b+1)^{1/2}\\
&\leq C(b+x)^{1/2}.
\end{align*}
This estimates separately the $L(x-y)$, $L(x+y)$, and two boundary
kernel terms $L(-y)$ and $L(y)$.

On $E_\sim$, the direct convolution term is estimated without losing
the diagonal cancellation:
\begin{align*}
\int_{x/2}^{2x}\lvert L(x-y)\rvert(b+y)^{1/2}\,\mathrm{d}y
&\leq
(b+x+1)^{1/2}
\int_{\mathbb R}\lvert L(z)\rvert\langle z\rangle^{1/2}\,\mathrm{d}z \leq C(b+x)^{1/2}.
\end{align*}
The remaining three terms are away from the origin of the Schwartz
kernel. Thus
\begin{align*} \int_{x/2}^{2x}
\bigl(\lvert L(x+y)\rvert+  \lvert L(-y)\rvert+ & \lvert L(y)\rvert\bigr)
(b+y)^{1/2}\,\mathrm{d}y\\
\leq & 
C_mx^{1-m}(b+x)^{1/2}
\leq C(b+x)^{1/2}.
\end{align*}
Consequently,
\begin{equation*}
J_\sim^-(x,b)+J_\sim^+(x,b)\leq C(b+x)^{1/2}.
\end{equation*}

Finally, if $y\in E_>$, then
$\lvert\theta x-y\rvert\geq y/2$ and
$\theta x+y\geq y$ for $0\leq\theta\leq1$. Applying the mean-value
theorem before integration and using the Schwartz decay of $L'$ gives
\begin{align*}
J_>^-(x,b)+J_>^+(x,b)
&\leq
C_mx\int_{2x}^\infty y^{-m}(b+y)^{1/2}\,\mathrm{d}y\\
&\leq C_m\left(b^{1/2}x^{2-m}+x^{5/2-m}\right) \leq C(b+x)^{1/2}.
\end{align*}
The last line holds because $x\geq1$ and $m>3$. We have therefore
proved, for every $x>0$,
\begin{equation*}
J_j^-(x,b)+J_j^+(x,b)\leq C(b+x)^{1/2},
\qquad j\in\{<,\sim,>\}.
\end{equation*}
Since
$\Lambda^{-1/2}(b+x)^{1/2}=(R+\rho)^{1/2}$, summing the six regional
bounds proves \eqref{eq:chainannularboundary}.

On odd functions, $\lvert D_\rho\rvert=\mathcal H\partial_\rho$ up to
a harmless sign. The Hilbert transform is skew-adjoint, and
$\widetilde\Pi_\Lambda\Gamma=\Gamma$. Therefore
\begin{align}
\left\lvert\int_0^\infty
\lvert D_\rho\rvert A(\rho)\Gamma(\rho)\,\mathrm{d}\rho\right\rvert
=\left\lvert\int_0^\infty
A'(\rho)T_\Lambda\Gamma(\rho)\,\mathrm{d}\rho\right\rvert =\left\lvert\int_0^\infty
A'(\rho)\bigl(T_\Lambda\Gamma(\rho)
-T_\Lambda\Gamma(0)\bigr)\,\mathrm{d}\rho\right\rvert.
\label{eq:chainpairing}
\end{align}
Here $\int_0^\infty A'(\rho)\,\mathrm{d}\rho=0$, because
$A(0)=0$ and $A$ has compact support. Since
$F(z)=\lvert z\rvert^{p-1}z$ and $H=rU$, differentiation gives
\begin{align*}
A'
&=F(U)+rF'(U)U'
=p\lvert U\rvert^{p-1}H'-(p-1)F(U),
\end{align*}
and
\begin{align*}
\lvert A'\rvert
&\leq C_p\lvert U\rvert^{p-1}
\left(\lvert H'\rvert+\frac{\lvert H\rvert}{r}\right).
\end{align*}
Moreover,
$
H'=\partial_\rho\lvert D_\rho\rvert^{-1}B
$
is, up to sign, the Hilbert transform of $B$. Applying the weighted
Hilbert-transform estimate proved in Step~1 to the odd extension of
$B$, followed by the one-sided Hardy estimate from Step~1 with
$g=\lvert H'\rvert$ and
$\lvert H(\rho)\rvert\leq\int_0^\rho\lvert H'(s)\rvert\,\mathrm{d}s$,
we have
\begin{equation}
\left\lVert r^{-\beta}H'\right\rVert_{L^a_\rho}
+\left\lVert r^{-\beta-1}H\right\rVert_{L^a_\rho}
\leq C_p\left\lVert r^{-\beta}B\right\rVert_{L^a_\rho}.
\label{eq:chainHBbound}
\end{equation}
Combining \eqref{eq:chainannularboundary},
\eqref{eq:chainpairing}, and \eqref{eq:chainHBbound}, and then applying
H\"older's inequality with the exponents in
\eqref{eq:chainHolderexponents}, we obtain
\begin{align*}
\left\lvert\int_0^\infty
\lvert D_\rho\rvert A(\rho)\Gamma(\rho)\,\mathrm{d}\rho\right\rvert
\leq &\;  C_p Z
\int_0^\infty r^{1/2}\lvert U\rvert^{p-1}
\left(\lvert H'\rvert+\frac{\lvert H\rvert}{r}\right)\,\mathrm{d}\rho\\
=&\; C_p Z
\int_0^\infty
\left(r^{(p-1)/q}\lvert U\rvert^{p-1}\right)
\left(r^{-\beta}\lvert H'\rvert
+r^{-\beta-1}\lvert H\rvert\right)\,\mathrm{d}\rho\\
\leq &\; C_p Z
\left\lVert r^{(p-1)/q}\lvert U\rvert^{p-1}
\right\rVert_{L^{a'}_\rho}
\left\lVert r^{-\beta}B\right\rVert_{L^a_\rho}.
\end{align*}
By 
\eqref{eq:chainhalflineidentities}, we proves
\eqref{eq:weightedchainkernelclaim}. 

We now take $U=u=P_{\leq N}^{(R)}v$ and $G=g_N$. The weighted
Littlewood--Paley estimate and \eqref{eq:morawetznorm} give
\begin{equation}
\left\lVert r^{-1/q}u\right\rVert_{L^q(\Omega_R)}
\leq C_p\left\lVert r^{-1/q}v\right\rVert_{L^q(\Omega_R)}
=C_p\mathcal M(v)^{1/q}.
\label{eq:chainfirstfactor}
\end{equation}
Lemma~\ref{lem:weightedSobolev} gives
\begin{equation*}
\left\lVert r^\beta\sqrt{-\Delta_{\Omega_R}}u
\right\rVert_{L^a}
\leq C_p
\left\lVert u\right\rVert_{\dot H_D^{1+\sigma}}.
\end{equation*}
If $\widehat v(\lambda)$ denotes the distorted Fourier transform and
$\chi_{\leq N}$ is the multiplier of $P_{\leq N}^{(R)}$, then
\begin{align*}
\left\lVert u\right\rVert_{\dot H_D^{1+\sigma}}^2
&=\int_0^\infty
\lambda^{2(1+\sigma)}
\lvert\chi_{\leq N}(\lambda)\widehat v(\lambda)\rvert^2
\,\mathrm{d}\lambda\\
&\leq C N^{2\sigma}
\int_0^\infty\lambda^2\lvert\widehat v(\lambda)\rvert^2
\,\mathrm{d}\lambda
=C N^{2\sigma}\left\lVert\nabla v\right\rVert_2^2.
\end{align*}
Thus
\begin{equation}
\left\lVert r^\beta\sqrt{-\Delta_{\Omega_R}}u
\right\rVert_{L^a}
\leq C_pN^\sigma\left\lVert\nabla v\right\rVert_2.
\label{eq:weightedchainsecondfactor}
\end{equation}
Substituting \eqref{eq:chainfirstfactor} and
\eqref{eq:weightedchainsecondfactor} into
\eqref{eq:weightedchainkernelclaim} , we get
\eqref{eq:weightedchain}.

\smallskip
\noindent\emph{Step 3: the localized radial Hardy estimate.}
Let
\begin{equation*}
h(\rho):=(\mathcal U_Rf)(\rho)=(R+\rho)f(R+\rho),\quad
k:=\lvert D_\rho\rvert h
=\mathcal U_R\sqrt{-\Delta_{\Omega_R}}f.
\end{equation*}
The boundary condition gives $h(0)=0$. With
$\Pi_L:=\mathcal U_RP_L^{(R)}\mathcal U_R^{-1}$,  we have
\begin{equation}
\mathcal U_R(\partial_rf)
=h'-\frac{h}{R+\rho},\qquad
\mathcal U_R(f/r)=\frac{h}{R+\rho}.
\label{eq:Hardyconjugatedderivatives}
\end{equation}
Write $k=\sum_Mk_M$, where $k_M:=\Pi_Mk$, and set
$h_M:=\lvert D_\rho\rvert^{-1}k_M$. It suffices to prove, for
$\nu=0,1$, that
\begin{equation}
\sup_{\rho>0}(R+\rho)^{-1/2}
\left\lvert T_{N,M}^{(\nu)}k_M(\rho)\right\rvert
\leq C\omega(N/M)
\sup_{\tau>0}(R+\tau)^{-1/2}\lvert k_M(\tau)\rvert,
\label{eq:Hardykernelclaim}
\end{equation}
where
\begin{align*}
T_{N,M}^{(0)}k_M
=\Pi_N\left(\frac{h_M}{R+\rho}\right),\qquad 
T_{N,M}^{(1)}k_M
:=\Pi_N\left(h_M'-\frac{h_M}{R+\rho}\right).
\end{align*}

We prove \eqref{eq:Hardykernelclaim} by scaling. Put
\begin{equation*}
s:=\frac NM,\quad b:=MR,\quad x:=M\rho,\quad y:=M\tau,
\end{equation*}
and define
\begin{equation*}
\varkappa(y):=k_M(y/M),\qquad
\widehat h(y):=Mh_M(y/M).
\end{equation*}
Then we have
\begin{equation*}
\lvert D_y\rvert\widehat h=\varkappa,\qquad
h_M'(\rho)=\widehat h'(x),\qquad
\frac{h_M(\rho)}{R+\rho}=\frac{\widehat h(x)}{b+x}.
\end{equation*}
The odd extension of $\varkappa$ has frequency support in a fixed unit
annulus. Define
\begin{equation*}
A_b:=\sup_{y>0}(b+y)^{-1/2}\lvert\varkappa(y)\rvert
\end{equation*}
and
\begin{equation*}
q_0(y):=\frac{\widehat h(y)}{b+y},\qquad
q_1(y):=\widehat h'(y)-\frac{\widehat h(y)}{b+y}.
\end{equation*}
If $\mathcal P_s$ denotes the half-line sine multiplier at frequency
$s$, the scaled form of \eqref{eq:Hardykernelclaim} is
\begin{equation}
\sup_{x>0}(b+x)^{-1/2}
\left\lvert\mathcal P_sq_\nu(x)\right\rvert
\leq C\omega(s)A_b,\qquad \nu=0,1.
\label{eq:Hardyscaledclaim}
\end{equation}

We now prove \eqref{eq:Hardyscaledclaim}. Since the multipliers
$\lvert\xi\rvert^{-1}$ and
$i\xi\lvert\xi\rvert^{-1}$ are smooth on the unit annulus, there are
Schwartz kernels $J_0,J_1$ such that
\begin{equation*}
\widehat h=J_0*\varkappa,\qquad
\widehat h'=J_1*\varkappa
\end{equation*}
on the odd extensions. For every Schwartz function $J$, we have
\begin{align*}
\int_{\mathbb R}\lvert J(y-z)\rvert(b+\lvert z\rvert)^{1/2}
\,\mathrm{d}z
 \leq
(b+\lvert y\rvert+1)^{1/2}
\int_{\mathbb R}\lvert J(y-z)\rvert
\langle y-z\rangle^{1/2}\,\mathrm{d}z
\leq C_J(b+\lvert y\rvert+1)^{1/2}.
\end{align*}
It follows that, for $y>0$,
\begin{equation}
\lvert\widehat h(y)\rvert+\lvert\widehat h'(y)\rvert
\leq C A_b(b+y+1)^{1/2}.
\label{eq:Hardyunitband}
\end{equation}
Oddness gives $\widehat h(0)=0$. Hence, when $0<y\leq1$, we have
\begin{equation*}
\lvert\widehat h(y)\rvert
\leq\int_0^y\lvert\widehat h'(z)\rvert\,\mathrm{d}z
\leq C A_by(b+y+1)^{1/2}.
\end{equation*}
Combining this with \eqref{eq:Hardyunitband} for $y\geq1$, we get
\begin{equation}
\lvert\widehat h(y)\rvert
\leq C A_b\min\{1,y\}(b+y+1)^{1/2}
\leq C A_b(b+y)^{1/2}.
\label{eq:Hardyboundarygain}
\end{equation}
Consequently,
\begin{equation}
\lvert q_0(y)\rvert+\lvert q_1(y)\rvert
\leq C A_b(b+y+1)^{1/2}.
\label{eq:Hardyqpointwise}
\end{equation}

Let
\begin{equation*}
Q_\nu(x):=\int_0^xq_\nu(y)\,\mathrm{d}y.
\end{equation*}
For $0<x\leq1$, \eqref{eq:Hardyunitband} and the boundary condition give
\begin{align*}
\lvert Q_0(x)\rvert
&\leq C A_b(b+1)^{1/2}
\int_0^x\frac{y}{b+y}\,\mathrm{d}y\\
&\leq C A_bx(b+1)^{1/2}
\leq C A_b(b+x)^{1/2},
\end{align*}
because $x^2(b+1)\leq b+x$. For $x\geq1$,
\eqref{eq:Hardyboundarygain} gives
\begin{align*}
\lvert Q_0(x)\rvert
&\leq\lvert Q_0(1)\rvert+
C A_b\int_1^x(b+y)^{-1/2}\,\mathrm{d}y\\
&\leq C A_b(b+1)^{1/2}
+2C A_b\bigl((b+x)^{1/2}-(b+1)^{1/2}\bigr)\\
&\leq C A_b(b+x)^{1/2}.
\end{align*}
Since
$Q_1(x)=\widehat h(x)-Q_0(x)$, the same bound and
\eqref{eq:Hardyboundarygain} imply
\begin{equation}
\lvert Q_0(x)\rvert+\lvert Q_1(x)\rvert
\leq C A_b(b+x)^{1/2},\qquad x>0.
\label{eq:Hardyprimitive}
\end{equation}

Let $\psi$ be the even Schwartz kernel associated with the smooth
annular sine projection. Then we have
\begin{equation*}
K_s(x,y):=
s\bigl[\psi(s(x-y))-\psi(s(x+y))\bigr],
\quad
\mathcal P_sq(x)=\int_0^\infty K_s(x,y)q(y)\,\mathrm{d}y.
\end{equation*}
For every $L>0$, the triangle inequality when $sx\geq1$, and the
mean-value theorem when $sx<1$, give
\begin{align}
\lvert K_s(x,y)\rvert
&\leq C_Ls\min\{1,sx\}
\bigl(\langle s(x-y)\rangle^{-L}
+\langle s(x+y)\rangle^{-L}\bigr),\notag\\
\lvert\partial_yK_s(x,y)\rvert
&\leq C_Ls^2\min\{1,sx\}
\bigl(\langle s(x-y)\rangle^{-L}
+\langle s(x+y)\rangle^{-L}\bigr).
\label{eq:Hardysinekernels}
\end{align}
For the first line in the case $sx<1$, we  uses the evenness of $\psi$
to write
\begin{equation*}
\psi(s(x-y))-\psi(s(x+y))
=\int_{-x}^{x}s\psi'(s(t-y))\,\mathrm{d}t.
\end{equation*}
For the second line, the oddness of $\psi'$ gives
\begin{equation*}
\psi'(s(x-y))+\psi'(s(x+y))
=\int_{-x}^{x}s\psi''(s(t-y))\,\mathrm{d}t.
\end{equation*}
These formulas prove the factor $sx$ in
\eqref{eq:Hardysinekernels}; the stated decay follows by increasing
$L$.

Taking $L>3$ in \eqref{eq:Hardysinekernels} and changing variables
$z=s(y-x)$ in the direct term gives
\begin{align}
\int_0^\infty\lvert K_s(x,y)\rvert(b+y+1)^{1/2}\,\mathrm{d}y
&\leq C\min\{1,sx\}
(b+x+1+s^{-1})^{1/2},\notag\\
\int_0^\infty\lvert\partial_yK_s(x,y)\rvert(b+y)^{1/2}\,\mathrm{d}y
&\leq Cs\min\{1,sx\}
(b+x+s^{-1})^{1/2}.
\label{eq:Hardykernelintegrals}
\end{align}
For example, the first inequality follows from
\begin{equation*}
(b+y+1)^{1/2}
\leq(b+x+1+s^{-1})^{1/2}
\langle s(x-y)\rangle^{1/2}
\end{equation*}
in the direct term. In the reflected term one uses the corresponding
inequality with $\langle s(x+y)\rangle^{1/2}$. The second inequality is
proved in the same way, with the additional factor $s$ from
$\partial_yK_s$.

Suppose first that $0<s\leq1$. Since $Q_\nu(0)=0$, the kernel is
Schwartz as $y\to\infty$, and \eqref{eq:Hardyprimitive} holds,
integration by parts gives
\begin{align*}
\mathcal P_sq_\nu(x)
&=-\int_0^\infty\partial_yK_s(x,y)Q_\nu(y)\,\mathrm{d}y,\\
\lvert\mathcal P_sq_\nu(x)\rvert
&\leq C A_bs\min\{1,sx\}
(b+x+s^{-1})^{1/2}.
\end{align*}
Let $m:=\min\{1,sx\}$. If $sx\geq1$, then
$b+x\geq s^{-1}$; if $sx<1$, then
$x^2/(b+x)\leq x$. Therefore
\begin{align*}
m^2\left(1+\frac{s^{-1}}{b+x}\right)
\leq & 
\begin{cases}
2,&sx\geq1,\\
s^2x^2+sx,&sx<1
\end{cases}\\
\leq & \; 2.
\end{align*}
It follows that
\begin{equation}
\sup_{x>0}(b+x)^{-1/2}
\lvert\mathcal P_sq_\nu(x)\rvert
\leq C sA_b,\quad 0<s\leq1.
\label{eq:Hardylowoutputbound}
\end{equation}

Suppose next that $s\geq1$. By
\eqref{eq:Hardyqpointwise} and the first inequality in
\eqref{eq:Hardykernelintegrals}, we have
\begin{equation*}
\lvert\mathcal P_sq_\nu(x)\rvert
\leq C A_bm(b+x+1+s^{-1})^{1/2}
\leq C A_bm(b+x+2)^{1/2}.
\end{equation*}
If $sx\geq1$, then $(b+x)^{-1}\leq s$; if $sx<1$, then
$x^2/(b+x)\leq x<s^{-1}$. Since $s\geq1$,
\begin{align*}
m^2\left(1+\frac2{b+x}\right)
\leq & 
\begin{cases}
1+2s,&sx\geq1,\\
s^2x^2+2s^2x,&sx<1
\end{cases}\\
\leq& \; 3s.
\end{align*}
Consequently,
\begin{equation}
\sup_{x>0}(b+x)^{-1/2}
\lvert\mathcal P_sq_\nu(x)\rvert
\leq C s^{1/2}A_b,\qquad s\geq1.
\label{eq:Hardyhighoutputbound}
\end{equation}
Equations \eqref{eq:Hardylowoutputbound} and
\eqref{eq:Hardyhighoutputbound} prove
\eqref{eq:Hardyscaledclaim}.

It remains to undo the scaling. One has
\begin{align*}
A_b
=M^{-1/2}
\sup_{\tau>0}(R+\tau)^{-1/2}\lvert k_M(\tau)\rvert,
\end{align*}
\begin{align*}
\sup_{\rho>0}(R+\rho)^{-1/2}
\lvert T_{N,M}^{(\nu)}k_M(\rho)\rvert
=M^{1/2}
\sup_{x>0}(b+x)^{-1/2}
\lvert\mathcal P_sq_\nu(x)\rvert.
\end{align*}
Thus \eqref{eq:Hardyscaledclaim} implies
\eqref{eq:Hardykernelclaim}. Moreover, for every radial $z$, we have
\begin{equation*}
\left\lVert r^{1/2}P_N^{(R)}z\right\rVert_{L_r^\infty}
=\sup_{\rho>0}(R+\rho)^{-1/2}
\left\lvert\Pi_N\mathcal U_Rz(\rho)\right\rvert.
\end{equation*}
Summing \eqref{eq:Hardykernelclaim} over $M$ and using
\eqref{eq:Hardyconjugatedderivatives},  we have
\eqref{eq:localizedHardy}.

\smallskip
\noindent\emph{Step 4: summation in the weighted Besov norm.}
Put $\alpha:=4/q$. Since $q\in(4,6)$,
\begin{equation*}
\frac23<\alpha<1,\qquad
1-\alpha>0,\qquad
\alpha-\frac12>0.
\end{equation*}
Multiplying \eqref{eq:localizedHardy} by $N^{-\alpha}$, summing first
over finite dyadic sets, and then using Tonelli's theorem, we have
\begin{align*}
\sum_NN^{-\alpha} & 
\sum_M\omega(N/M)
\left\lVert r^{1/2}P_M^{(R)}
\sqrt{-\Delta_{\Omega_R}}f\right\rVert_{L_r^\infty}\\
= & 
\sum_MM^{-\alpha}
\left\lVert r^{1/2}P_M^{(R)}
\sqrt{-\Delta_{\Omega_R}}f\right\rVert_{L_r^\infty} \times
\left[
\sum_{N\leq M}\left(\frac NM\right)^{1-\alpha}
+\sum_{N>M}\left(\frac NM\right)^{1/2-\alpha}
\right].
\end{align*}
The bracket is bounded uniformly in $M$, because
\begin{equation*}
\sum_{j\leq0}2^{(1-\alpha)j}<\infty,\qquad
\sum_{j\geq1}2^{-(\alpha-1/2)j}<\infty.
\end{equation*}
Passing to the monotone limit, we have
\begin{equation*}
B_{\partial_rf}+B_{f/r}
\leq C
\sum_MM^{-4/q}
\left\lVert r^{1/2}P_M^{(R)}
\sqrt{-\Delta_{\Omega_R}}f\right\rVert_{L_r^\infty}
=C B_{\sqrt{-\Delta_{\Omega_R}}f},
\end{equation*}
which is \eqref{eq:radialHardyBproved}, and completes the proof.
\end{proof}

For a spacetime function $z(t,x)$, put $q=p+1\in (4,6)$, we define the weighted Strichartz norm
\begin{equation}
 \mathcal{B}_p(z):=
 \sum_{N\in2^\mathbb Z}N^{-4/q}
 \left\lVert r^{1/2}P_N^{(R)}(\partial_t,\sqrt{-\Delta_{\Omega_R}})z\right\rVert
 _{L_t^{q/2}L_r^\infty}.
 \label{eq:Bpnorm}
\end{equation}

\begin{lemma}
\label{lem:forcedBesov}
Suppose
\begin{equation*}
\left\{
\begin{aligned}
& \partial_t^2y-\Delta y=G,
&&(t,x)\in I\times\Omega_R,\\
&y(t,x)=0,&&x\in\partial\Omega_R,\\
&(y,\partial_ty)|_{t=0}=(y_0,y_1),
\end{aligned}
\right.
\end{equation*}
Then we have
\begin{align}
 \mathcal{B}_p(y)\lesssim_p{}&
 \left\lVert P_{\ge1}(y_0,y_1)\right\rVert
 _{\dot H_D^{s_c}\times\dot H_D^{s_c-1}}
  +\left\lVert P_{\ge1}G\right\rVert_{L_t^1\dot H_D^{s_c-1}}
 \notag\\
 &\qquad +\left\lVert P_{\le1}(y_0,y_1)\right\rVert_{L^2\times\dot H_D^{-1}}
 +\left\lVert P_{\le1}G\right\rVert_{L_t^1\dot H_D^{-1}}.
 \label{eq:forcedBesovEstimate}
\end{align}
\end{lemma}

\begin{proof}
We separate the argument into the estimate at one dyadic frequency and
the summation of those estimates.  Since $q=p+1\in(4,6)$, the time
exponent
\begin{equation*}
 b=\frac q2=\frac{p+1}{2}
\end{equation*}
satisfies $2<b<3$.  If $f_N=P_N^{(R)}f$, then, by Euler's formula, Lemma~\ref{lem:weightedStrichartz} gives
\begin{equation}
 \left\lVert r^{1/2}\cos(t\sqrt{-\Delta_{\Omega_R}})f_N\right\rVert_{L_t^{q/2}L_r^\infty}
 +\left\lVert r^{1/2}\sin(t\sqrt{-\Delta_{\Omega_R}})f_N\right\rVert_{L_t^{q/2}L_r^\infty}
 \lesssim_p N^{1-2/q}\left\lVert f_N\right\rVert_2.
 \label{eq:cosinesineweighted}
\end{equation}
The estimate is unchanged when $t$ is replaced by $t-s$, since the
$L_t^{q/2}$ norm is invariant under time translation.

The frequency projections commute with $\sqrt{-\Delta_{\Omega_R}}$ and with the wave group.
Consequently, the dyadic Duhamel formula is
\begin{align*}
 P_Ny(t)
 =&\ \cos(t\sqrt{-\Delta_{\Omega_R}})P_Ny_0+\frac{\sin(t\sqrt{-\Delta_{\Omega_R}})}{\sqrt{-\Delta_{\Omega_R}}} 
 P_Ny_1 +\int_0^t\frac{\sin((t-s)\sqrt{-\Delta_{\Omega_R}})}{\sqrt{-\Delta_{\Omega_R}}} 
 P_NG(s)\,\mathrm{d} s
\end{align*}
and differentiation gives
\begin{align}
 \partial_tP_Ny(t)
 =&-\sqrt{-\Delta_{\Omega_R}}\sin(t\sqrt{-\Delta_{\Omega_R}})P_Ny_0+\cos(t\sqrt{-\Delta_{\Omega_R}})P_Ny_1\notag\\
 &\quad +\int_0^t\cos((t-s)\sqrt{-\Delta_{\Omega_R}})P_NG(s)\,\mathrm{d} s,
 \label{eq:forcedDuhameltime}\\
 \sqrt{-\Delta_{\Omega_R}} P_Ny(t)
 =&\ \sqrt{-\Delta_{\Omega_R}}\cos(t\sqrt{-\Delta_{\Omega_R}})P_Ny_0+\sin(t\sqrt{-\Delta_{\Omega_R}})P_Ny_1\notag\\
 &\quad +\int_0^t\sin((t-s)\sqrt{-\Delta_{\Omega_R}})P_NG(s)\,\mathrm{d} s.
 \label{eq:forcedDuhamelspace}
\end{align}
On the support of $P_N^{(R)}$ one has
$\left\lVert \sqrt{-\Delta_{\Omega_R}} P_Nf\right\rVert_2\lesssim N\left\lVert P_Nf\right\rVert_2$.  Hence
\eqref{eq:cosinesineweighted} applied to the homogeneous terms in
\eqref{eq:forcedDuhameltime}--\eqref{eq:forcedDuhamelspace} yields
\begin{align}
 &\left\lVert r^{1/2}P_N(\partial_t,\sqrt{-\Delta_{\Omega_R}})y_{\mathrm{hom}}\right\rVert_
 {L_t^{q/2}L_r^\infty} \lesssim_p
 N^{1-2/q}\left(
 N\left\lVert P_Ny_0\right\rVert_2+\left\lVert P_Ny_1\right\rVert_2\right).
 \label{eq:homogeneousdyadicBesov}
\end{align}

We next justify the retarded Duhamel term without changing the order of
the mixed norm.  On $t\ge0$, Minkowski's integral inequality and
\eqref{eq:cosinesineweighted} give
\begin{align*}
 &\left\lVert r^{1/2}\int_0^t
  (\cos((t-s)\sqrt{-\Delta_{\Omega_R}}),\sin((t-s)\sqrt{-\Delta_{\Omega_R}}))P_NG(s)\,\mathrm{d} s\right\rVert_
  {L_t^{q/2}([0,\infty))L_r^\infty}\notag\\
 &\quad\le
 \int_0^\infty
 \left\lVert \mathbf 1_{\{t\ge s\}}r^{1/2}
 (\cos((t-s)\sqrt{-\Delta_{\Omega_R}}),\sin((t-s)\sqrt{-\Delta_{\Omega_R}}))P_NG(s)\right\rVert_
 {L_t^{q/2}L_r^\infty}\,\mathrm{d} s\notag\\
 &\quad\lesssim_p
 N^{1-2/q}\int_0^\infty\left\lVert P_NG(s)\right\rVert_2\,\mathrm{d} s.
\end{align*}
The indicator only restricts the time integral, so its removal can only
increase the norm.  For $t\le0$, write
$\int_0^t=-\int_t^0$ and repeat the argument with
$\mathbf 1_{\{t\le s\}}$ and $s\le0$.  Combining the two time half-lines
gives
\begin{align}
 &\left\lVert r^{1/2}P_N(\partial_t,\sqrt{-\Delta_{\Omega_R}})y_{\mathrm{Duh}}\right\rVert_
 {L_t^{q/2}L_r^\infty}
 \lesssim_p
 N^{1-2/q}\left\lVert P_NG\right\rVert_{L_t^1L_x^2}.
 \label{eq:retardedDuhamelMinkowski}
\end{align}

Adding \eqref{eq:homogeneousdyadicBesov} and
\eqref{eq:retardedDuhamelMinkowski}, and then multiplying by the weight
$N^{-4/q}$ from \eqref{eq:Bpnorm}, gives the fundamental dyadic estimate
\begin{align}
 &N^{-4/q}\left\lVert r^{1/2}P_N(\partial_t,\sqrt{-\Delta_{\Omega_R}})y\right\rVert_
 {L_t^{q/2}L_r^\infty}\notag\\
 &\qquad\lesssim
 N^{2-6/q}\left\lVert P_Ny_0\right\rVert_2
 +N^{1-6/q}\left\lVert P_Ny_1\right\rVert_2
 +N^{1-6/q}\left\lVert P_NG\right\rVert_{L_t^1L_x^2}.
 \label{eq:dyadicforcedbound}
\end{align}

We first sum \eqref{eq:dyadicforcedbound} over the high frequencies
$N\ge1$.  Factoring out the critical Sobolev weights, its right-hand
side becomes
\begin{align*}
 N^{2-s_c-6/q}\bigl(
 N^{s_c}\left\lVert P_Ny_0\right\rVert_2
 +N^{s_c-1}\left\lVert P_Ny_1\right\rVert_2
 +N^{s_c-1}\left\lVert P_NG\right\rVert_{L_t^1L_x^2}\bigr).
\end{align*}
The exponent left outside the parentheses is strictly negative.  More
precisely,
\begin{equation}\label{def:etap}
 2-s_c-\frac6q=-\eta_p,
 \qquad
 \eta_p:=s_c+\frac6q-2
 =\frac{(p-3)(5-p)}{2(p^2-1)}>0.
\end{equation}
Both strict inequalities $p>3$ and $p<5$ are used here.  Since the
frequencies are dyadic, $\eta_p>0$ implies
$\sum_{N\ge1}N^{-2\eta_p}<\infty$.  The Cauchy--Schwarz inequality therefore gives that
\begin{align}
 \sum_{N\geq1}N^{-\eta_p}N^{s_c}\left\lVert P_Ny_0\right\rVert_2
 &\leq
 \left(\sum_{N\geq1}N^{-2\eta_p}\right)^{1/2}
 \left(\sum_{N\geq1}N^{2s_c}\left\lVert P_Ny_0\right\rVert_2^2\right)^{1/2}\lesssim_p\left\lVert P_{\geq1}y_0\right\rVert_{\dot H_D^{s_c}}.
 \label{eq:highfrequencysum}
\end{align}
The same calculation with $N^{s_c-1}\left\lVert P_Ny_1\right\rVert_2$ gives
\begin{equation}
 \sum_{N\ge1}N^{-\eta_p}N^{s_c-1}\left\lVert P_Ny_1\right\rVert_2
 \lesssim_p
 \left\lVert P_{\ge1}y_1\right\rVert_{\dot H_D^{s_c-1}}.
 \label{eq:highvelocitysum}
\end{equation}
For the forcing, we  interchange the dyadic square sum and the
$L_t^1$ integral.  The vector-valued Minkowski inequality gives
\begin{align}
 &\left(\sum_{N\ge1}N^{2(s_c-1)}
 \left\lVert P_NG\right\rVert_{L_t^1L_x^2}^2\right)^{1/2}  \le
 \int_\mathbb{R}\left(\sum_{N\ge1}N^{2(s_c-1)}
 \left\lVert P_NG(t)\right\rVert_2^2\right)^{1/2}\,\mathrm{d} t \lesssim
 \left\lVert P_{\ge1}G\right\rVert_{L_t^1\dot H_D^{s_c-1}}.
 \label{eq:highforcingMinkowski}
\end{align}
Combining this with the Cauchy--Schwarz inequality in $N$ yields
\begin{equation}
 \sum_{N\ge1}N^{-\eta_p}N^{s_c-1}
 \left\lVert P_NG\right\rVert_{L_t^1L_x^2}
 \lesssim_p
 \left\lVert P_{\ge1}G\right\rVert_{L_t^1\dot H_D^{s_c-1}}.
 \label{eq:highforcingsum}
\end{equation}

It remains to sum the low frequencies.  For $N\le1$, rewrite
\eqref{eq:dyadicforcedbound} as

\begin{align}
 N^{\mu_p}\bigl(
 \left\lVert P_Ny_0\right\rVert_2+N^{-1}\left\lVert P_Ny_1\right\rVert_2
 +N^{-1}\left\lVert P_NG\right\rVert_{L_t^1L_x^2}\bigr).
 \label{eq:lowforcedbound}
\end{align}
Here
\begin{equation*}
 \mu_p:=2-\frac6q=\frac{2(p-2)}{p+1}>0.
\end{equation*}
Thus $\sum_{N\le1}N^{2\mu_p}<\infty$.  Applying the Cauchy--Schwarz inequality to
\eqref{eq:lowforcedbound} and then using the frequency square-function
identity \eqref{eq:squarefunctionL2} gives
\begin{align}
 \sum_{N\leq1}N^{2-6/q}\left\lVert P_Ny_0\right\rVert_2
 &\lesssim_p\left\lVert P_{\leq1}y_0\right\rVert_2,\label{eq:lowpositionsum}\\
 \sum_{N\leq1}N^{1-6/q}\left\lVert P_Ny_1\right\rVert_2
 &\lesssim_p\left\lVert P_{\leq1}y_1\right\rVert_{\dot H_D^{-1}}.
 \label{eq:lowvelocitysum}
\end{align}
For clarity, the forcing term uses the same vector-valued Minkowski step
as in \eqref{eq:highforcingMinkowski}, now at regularity $-1$:
\begin{align*}
 &\left(\sum_{N\le1}N^{-2}
 \left\lVert P_NG\right\rVert_{L_t^1L_x^2}^2\right)^{1/2}\notag \le
 \int_\mathbb{R}\left(\sum_{N\le1}N^{-2}
 \left\lVert P_NG(t)\right\rVert_2^2\right)^{1/2}\,\mathrm{d} t
 \lesssim
 \left\lVert P_{\le1}G\right\rVert_{L_t^1\dot H_D^{-1}}.
\end{align*}
Consequently,
\begin{equation}
 \sum_{N\le1}N^{\mu_p-1}\left\lVert P_NG\right\rVert_{L_t^1L_x^2}
 \lesssim_p
 \left\lVert P_{\le1}G\right\rVert_{L_t^1\dot H_D^{-1}}.
 \label{eq:lowforcingsum}
\end{equation}

Finally, summing the left-hand side of \eqref{eq:dyadicforcedbound} over
all dyadic $N$ is $\mathcal{B}_p(y)$.  The high-frequency estimates
\eqref{eq:highfrequencysum}, \eqref{eq:highvelocitysum}, and
\eqref{eq:highforcingsum}, together with the low-frequency estimates
\eqref{eq:lowpositionsum}, \eqref{eq:lowvelocitysum}, and
\eqref{eq:lowforcingsum}, prove \eqref{eq:forcedBesovEstimate}.
\end{proof}

\begin{proposition}
\label{prop:smallBesov}
Suppose that $w$ is the small global solution of \eqref{eq:nlw} and its
initial data are frequency supported in $\{\lambda\ge1\}$. If
\begin{equation*}
 \left\lVert (w_0,w_1)\right\rVert_{\dot H_D^{s_c}\times\dot H_D^{s_c-1}}
 \le\varepsilon\le\varepsilon_p,
\end{equation*}
then
\begin{equation}
 \mathcal{B}_p(w)
 +\left\lVert w\right\rVert_{L_t^\infty L_x^2}
 +\left\lVert \partial_tw\right\rVert_{L_t^\infty\dot H_D^{-1}}
 \lesssim_p\varepsilon.
 \label{eq:smallBesov}
\end{equation}
\end{proposition}

\begin{proof}
The frequency support assumption gives
\begin{equation*}
 \left\lVert w_0\right\rVert_2+\left\lVert w_1\right\rVert_{\dot H_D^{-1}}
 \lesssim\varepsilon.
\end{equation*}
Moreover,
\begin{align}
 \left\lVert F(w)\right\rVert_{L_t^1\dot H_D^{-1}}
 \lesssim
 \left\lVert F(w)\right\rVert_{L_t^1L_x^{6/5}}\lesssim
 \left\lVert w\right\rVert_{L_t^\infty L_x^2}
 \left\lVert w\right\rVert_{L_t^{p-1}L_x^{3(p-1)}}^{p-1}.
 \label{eq:Hminusoneforce}
\end{align}
Plancherel applied to the Duhamel formula at regularity zero gives, with
\begin{equation*}
 X:=\left\lVert w\right\rVert_{L_t^\infty L_x^2}
 +\left\lVert \partial_tw\right\rVert_{L_t^\infty\dot H_D^{-1}},
\end{equation*}
the estimate
\begin{equation*}
 X\lesssim\varepsilon+\left\lVert F(w)\right\rVert_{L_t^1\dot H_D^{-1}}
 \lesssim\varepsilon+\left\lVert w\right\rVert_{S_c}^{p-1}X.
\end{equation*}
By \eqref{eq:smalldata} and the choice of $\varepsilon_p$, the last term
is absorbed. This proves the last two bounds in \eqref{eq:smallBesov}.

Apply Lemma~\ref{lem:forcedBesov} to
$G=-F(w)$. The high-frequency forcing norm is bounded by
\eqref{eq:chain}, while the low-frequency norm is bounded by
\eqref{eq:Hminusoneforce}.  Hence
\begin{equation*}
 \mathcal{B}_p(w)\lesssim_p\varepsilon+C_p\varepsilon^p
 \lesssim_p\varepsilon.
\end{equation*}
Together with the already established $L^2\times\dot H_D^{-1}$ bound,
this proves \eqref{eq:smallBesov}.
\end{proof}

\section{Bourgain's high-low frequency truncation decomposition}
\label{sec:decomposition}

After the scale normalization of Section~\ref{sec:mainproof}, the frequency truncation
decomposition takes the form
\begin{equation*}
 (v_0,v_1):=P_{\le1}^{(R)}(u_0,u_1),
 \qquad
 (w_0,w_1):=P_{>1}^{(R)}(u_0,u_1),
\end{equation*}
where the high frequency part has critical norm at most $\varepsilon$ if $R$ is sufficiently large. Let $w$ be
the global small solution
\begin{equation*}
\left\{
\begin{aligned}
& \partial_t^2w-\Delta w+F(w)=0,
&&(t,x)\in I\times\Omega_R,\\
&w(t,x)=0,&&x\in\partial\Omega_R,\\
&(w,\partial_tw)|_{t=0}=(w_0,w_1),
\end{aligned}
\right.
\qquad
\label{eq:nlw w}
\end{equation*}
Write $u=v+w$. Then $v$ solves the difference equation
\begin{equation}
\left\{
\begin{aligned}
& \partial_t^2v-\Delta v+F(v)+\mathcal R(v,w)=0,
&&(t,x)\in I\times\Omega_R,\\
&v(t,x)=0,&&x\in\partial\Omega_R,\\
&(v,\partial_tv)|_{t=0}=(v_0,v_1)),
\end{aligned}
\right.
\qquad
 \label{eq:vequation}
\end{equation}
where $ \mathcal R(v,w):=F(v+w)-F(v)-F(w).$

By the elementary calculations, we have
\begin{align}
 \left\lvert \mathcal R(v,w)\right\rvert
 &\lesssim_p \left\lvert v\right\rvert^{p-1}\left\lvert w\right\rvert
                 +\left\lvert v\right\rvert\left\lvert w\right\rvert^{p-1},
 \label{eq:Taylor1}\\
 \left\lvert \mathcal R(v,w)-p\left\lvert v\right\rvert^{p-1}w\right\rvert
 &\lesssim_p \left\lvert v\right\rvert^{p-2}\left\lvert w\right\rvert^2
                 +\left\lvert v\right\rvert\left\lvert w\right\rvert^{p-1}.
 \label{eq:Taylor2}
\end{align}
Since the second-order remainder is essential for the corrected energy, we
include the proof here. In deed, the identities
\begin{equation*}
 F'(z)=p|z|^{p-1},\qquad
 |F'(z_1)-F'(z_2)|
 \lesssim_p(|z_1|+|z_2|)^{p-2}|z_1-z_2|,
\end{equation*}
together with the fundamental theorem of calculus give
\begin{align*}
 F(v+w)-F(v)
 &=w\int_0^1F'(v+\theta w)\,\mathrm{d}\theta,\\
 F(v+w)-F(v)-F'(v)w
 &=w\int_0^1\bigl(F'(v+\theta w)-F'(v)\bigr)\,\mathrm{d}\theta.
\end{align*}
The first identity, after subtracting $F(w)$ and separating the regions
$|w|\leq|v|$ and $|w|>|v|$, yields \eqref{eq:Taylor1}. The second gives
\begin{equation*}
 |F(v+w)-F(v)-p|v|^{p-1}w|
 \lesssim_p |v|^{p-2}|w|^2+|w|^p.
\end{equation*}
Subtracting $F(w)$ from the left-hand side and using
$|w|^p\mathbf1_{\{|w|\leq|v|\}}\leq|v||w|^{p-1}$ and the symmetric
estimate on $\{|w|>|v|\}$ proves \eqref{eq:Taylor2}.

For an energy solution $v$, we define
\begin{equation*}
 E(v)(t):=\frac12\left\lVert \partial_tv(t)\right\rVert_2^2
 +\frac12\left\lVert \nabla v(t)\right\rVert_2^2
 +\frac1{p+1}\left\lVert v(t)\right\rVert_{p+1}^{p+1},
\end{equation*}
which is not conservative since $v$ satisfies \eqref{eq:vequation}  with a forcing term $\mathcal R(v,w)$.

\begin{lemma}\label{lem:energyperturbation}
Let
\begin{equation*}
 Q_1:=\frac{2p}{p-3},
 \qquad b_1:=\frac{2p(p-1)}{p+3},
 \qquad \vartheta:=\frac{p+3}{2p}.
\end{equation*}
On every compact interval $J$, the small solution $w$ satisfies
\begin{equation}
 \left\lVert w\right\rVert_{L_t^{b_1}L_x^{2p}(J)}
 \lesssim
 \left\lVert w\right\rVert_{L_t^\infty\dot H_D^{s_c}(J)}^{1-\vartheta}
 \left\lVert w\right\rVert_{L_t^{p-1}L_x^{3(p-1)}(J)}^\vartheta.
 \label{eq:wenergyspace}
\end{equation}
If $v\in L_t^{Q_1}L_x^{2p}(J)$, then we have
\begin{align}
 \left\lVert F(v)\right\rVert_{L_t^1L_x^2(J)}
 &\le \left\lvert J\right\rvert^{1-p/Q_1}
 \left\lVert v\right\rVert_{L_t^{Q_1}L_x^{2p}(J)}^p,
 \label{eq:energyF}\\
 \left\lVert \mathcal R(v,w)\right\rVert_{L_t^1L_x^2(J)}
 &\lesssim_p
 \left\lvert J\right\rvert^{\kappa_p}
 \left\lVert v\right\rVert_{L_t^{Q_1}L_x^{2p}}^{p-1}
 \left\lVert w\right\rVert_{L_t^{b_1}L_x^{2p}} +
 \left\lVert v\right\rVert_{L_t^{Q_1}L_x^{2p}}
 \left\lVert w\right\rVert_{L_t^{b_1}L_x^{2p}}^{p-1},
 \label{eq:energyR}
\end{align}
where
$\kappa_p=1-(p-1)/Q_1-1/b_1>0$ for $p<5$.
\end{lemma}

\begin{proof}
By the Sobolev estimate, we have 
$\dot H_D^{s_c}\hookrightarrow L^{3(p-1)/2}$. Interpolation with the
$L_t^{p-1}L_x^{3(p-1)}$ component of $S_c$ gives
\eqref{eq:wenergyspace}, because
\begin{equation*}
 \frac1{2p}=\frac{1-\vartheta}{3(p-1)/2}
 +\frac{\vartheta}{3(p-1)},
 \qquad \frac1{b_1}=\frac{\vartheta}{p-1}.
\end{equation*}
Equation \eqref{eq:energyF} is H\"older's inequality in time. The Taylor
bound \eqref{eq:Taylor1}, the H\"older inequality with all factors in
$L^{2p}$, and H\"older's inequality in time give \eqref{eq:energyR}. 
\end{proof}

\begin{proposition} \label{prop:affine}
Assume $w$ satisfies Proposition~\ref{prop:smallBesov} and
$(v_0,v_1)\in(\dot H_D^1\cap L^{p+1})\times L^2$. There is a unique local
solution of \eqref{eq:vequation} in the energy space $\dot H_D^1\times L^2$. It continues as long as $E(v)(t)$ remains finite.
\end{proposition}

\begin{proof} It is corresponding to Proposition $4.2$ in \cite{XuYang2024} and Proposition $3.2$ in \cite{XuYangYou2025}. We include the proof for completeness. 
The pair $(Q_1,2p)$ satisfies \eqref{eq:radialadmissible} with $s=1$.
For $t_0\in J$, write
\begin{align*}
 \Psi(v)(t)
 =&\ \mathcal W_R(t-t_0)(v(t_0),\partial_tv(t_0)) -\int_{t_0}^t \frac{\sin((t-s)\sqrt{-\Delta_{\Omega_R}})}{\sqrt{-\Delta_{\Omega_R}}} 
 \bigl(F(v)+\mathcal R(v,w)\bigr)(s)\,\mathrm{d} s.
\end{align*}

For an energy-class function $v$ on $J$, we define 
\begin{equation*}
 \left\lVert v\right\rVert_{Y(J)}
 :=\left\lVert v\right\rVert_{L_t^\infty\dot H_D^1(J)}
  +\left\lVert \partial_tv\right\rVert_{L_t^\infty L_x^2(J)}
  +\left\lVert v\right\rVert_{L_t^{Q_1}L_x^{2p}(J)}.
\end{equation*}

Let
\begin{equation*}
 M:=1+\left\lVert \nabla v(t_0)\right\rVert_2+\left\lVert \partial_tv(t_0)\right\rVert_2.
\end{equation*}
After enlarging the implicit constant,
Lemma~\ref{lem:radialStrichartz} and Lemma~\ref{lem:energyperturbation} imply
\begin{align*}
 \left\lVert \Psi(v)\right\rVert_{Y(J)}
 \leq C M
 +C_p|J|^{1-p/Q_1}\left\lVert v\right\rVert_{L_t^{Q_1}L_x^{2p}}^p
 +C_p\mathcal R_J(v,w),
\end{align*}
where $\mathcal R_J(v,w)$ denotes the right-hand side of
\eqref{eq:energyR}. The same estimates applied to
$F(v)-F(\widetilde v)$ and
$\mathcal R(v,w)-\mathcal R(\widetilde v,w)$ yield
\begin{align*}
 \left\lVert \Psi(v)-\Psi(\widetilde v)\right\rVert_{Y(J)}
 \leq o_{|J|\to0}(1)\left\lVert v-\widetilde v\right\rVert_{Y(J)}
\end{align*}
uniformly for $v,\widetilde v$ in the radius-$M$ ball. Here the small
factor follows from $1-p/Q_1>0$, $\kappa_p>0$, and the absolute
continuity of $\left\lVert w\right\rVert_{L_t^{b_1}L_x^{2p}(J)}$. Thus $\Psi$ is a
contraction after shrinking $J$.

Smooth approximation justifies
\begin{equation*}
 E(v)(t)=E(v)(t_0)-\int_{t_0}^t
 \left\langle \partial_sv(s),\mathcal R(v,w)(s)\right\rangle\,\mathrm{d} s,
\end{equation*}
and passage to the limit follows from the $Y(J)$ convergence and
\eqref{eq:energyR}. If $E(v)$
stays finite, the $\dot H_D^1\times L^2$ norm and the potential-energy
norm stay finite, while \eqref{eq:wenergyspace} is finite on compact
intervals. The same contraction can therefore be restarted at the end of
$J$, proving the continuation criterion.
\end{proof}

\section{Modified energy and exterior Morawetz estimate}
\label{sec:modified}

In order to control the growth of the energy for $v$, we should modify the energy $E(v)(t)$ with the correction terms, one of which is the exterior Morawetz term.

\subsection{The exterior Morawetz identity}

Define
\begin{equation}
 M_R(v)(t):=
 \int_{\Omega_R}\partial_tv(t,x)
 \left(\partial_rv(t,x)+\frac{v(t,x)}r\right)\,\mathrm{d} x.
 \label{eq:Morawetzfunctional}
\end{equation}
Lemma~\ref{lem:basictools} gives
\begin{equation}
 \left\lvert M_R(v)(t)\right\rvert\lesssim E(v)(t).
 \label{eq:Mbound}
\end{equation}
\begin{lemma}
\label{lem:Morawetzidentity}
For a smooth radial Dirichlet solution of \eqref{eq:vequation}, we have
\begin{align}
 \frac{\,\mathrm{d}}{\,\mathrm{d} t}M_R(v)
 =&-2\pi R^2\left\lvert \partial_rv(t,R)\right\rvert^2
 -\frac{p-1}{p+1}
 \int_{\Omega_R}\frac{\left\lvert v\right\rvert^{p+1}}r\,\mathrm{d} x -\int_{\Omega_R}\mathcal R(v,w)
 \left(\partial_rv+\frac vr\right)\,\mathrm{d} x.
 \label{eq:Morawetzidentity}
\end{align}
and
\begin{equation}
 \frac{\,\mathrm{d}}{\,\mathrm{d} t}E(v)(t)
 =-\int_{\Omega_R}\partial_tv\,\mathcal R(v,w)\,\mathrm{d} x.
 \label{eq:energyvariation}
\end{equation}
The boundary term in \eqref{eq:Morawetzidentity} has a favorable sign.
For the energy solutions the identities hold in the integrated
inequality sense after smooth approximation.
\end{lemma}

\begin{proof}
Write $v=v(t,r)$ and $\,\mathrm{d} x=4\pi r^2\,\mathrm{d} r$. Substitute
$v_{tt}=v_{rr}+2r^{-1}v_r-F(v)-\mathcal R(v,w)$ into the derivative of
\eqref{eq:Morawetzfunctional}. The linear terms are
\begin{align*}
4\pi\int_R^\infty
 &(v_{rr}+2r^{-1}v_r)(v_r+v/r)r^2
 {}+v_t(v_{tr}+v_t/r)r^2\,\mathrm{d} r =-2\pi R^2\left\lvert v_r(t,R)\right\rvert^2.
\end{align*}
Here $v(t,R)=v_t(t,R)=0$, and all terms at infinity vanish by
approximation. For the pure nonlinearity,
\begin{align*}
 -4\pi\int_R^\infty F(v)(v_r+v/r)r^2\,\mathrm{d} r
 &=
 \frac{8\pi}{p+1}\int_R^\infty\left\lvert v\right\rvert^{p+1}r\,\mathrm{d} r
 -4\pi\int_R^\infty\left\lvert v\right\rvert^{p+1}r\,\mathrm{d} r\\
 &=-\frac{p-1}{p+1}
 \int_{\Omega_R}\frac{\left\lvert v\right\rvert^{p+1}}r\,\mathrm{d} x.
\end{align*}
The perturbation gives the last term in
\eqref{eq:Morawetzidentity}. Multiplication of \eqref{eq:vequation} by
$v_t$ gives \eqref{eq:energyvariation}. Smooth approximation and Fatou's
lemma preserve the favorable boundary sign for the energy solutions.
\end{proof}

Set
\begin{equation*}
 \mathcal{M}(v)(t):=\int_{\Omega_R}\frac{\left\lvert v(t,x)\right\rvert^{p+1}}r\,\mathrm{d} x.
\end{equation*}
Define the energy-Morawets quantity with another important correction term
\begin{equation}
 \mathcal{E}(v,w)(t):=
 E(v)(t)+cM_R(v)(t)
 +\int_{\Omega_R}F(v(t,x))w(t,x)\,\mathrm{d} x.
 \label{eq:modifiedenergy}
\end{equation}
where $c>0$ is a sufficiently small. 
The sign of the last term is chosen according to
\eqref{eq:energyvariation}. Its derivative cancels
\begin{equation*}
-p\int_{\Omega_R}\partial_tv\left\lvert v\right\rvert^{p-1}w\,\mathrm{d}x.
\end{equation*}

\begin{lemma}\label{lem:coercivity}
If $\left\lVert w\right\rVert_{L_t^\infty\dot H_D^{s_c}}\le\varepsilon$ and
$\varepsilon+c$ is sufficiently small depending only on $p$, then we have
\begin{equation*}
 \frac12E(v)(t)\le \mathcal{E}(v,w)(t)\le2E(v)(t).
\end{equation*}
\end{lemma}

\begin{proof}
By \eqref{eq:Mbound}, it remains to control the correction. Set
$b=2/(p-1)$ and $a=p-b$. Then
\begin{equation*}
 \frac{a}{p+1}+\frac b2=1,
 \qquad
 \frac1{3(p-1)/2}+\frac{a}{p+1}+\frac b6=1.
\end{equation*}
Hence, by the Sobolev inequality and Lemma~\ref{lem:radialsobolev}, we have
\begin{align*}
 \left|\int F(v)w\,\mathrm{d} x\right|
 &\le
 \left\lVert w\right\rVert_{L^{3(p-1)/2}}
 \left\lVert v\right\rVert_{L^{p+1}}^a\left\lVert v\right\rVert_{L^6}^b
  \lesssim_p
 \left\lVert w\right\rVert_{\dot H_D^{s_c}}E(v).
\end{align*}
This proves the claim.
\end{proof}

\subsection{The weighted paraproduct estimate}

\begin{lemma}\label{lem:paraproduct}
Let $v$ be radial with finite energy and let $z$ be radial and smooth.
Then we have
\begin{equation}
 \left\lvert \left\langle F(v),z\right\rangle\right\rvert
 \lesssim_p
 E(v)^{1/2}\mathcal{M}(v)^{\frac{p-1}{p+1}}B_z.
 \label{eq:paraproduct}
\end{equation}
where $B_z$ is defined by \eqref{def:Bz}.
\end{lemma}

\begin{proof}
Decompose $z=\sum_NP_N^{(R)}z$ and, for each $N$, write
\begin{equation*}
 F(v)=\bigl(F(v)-F(P_{\le N}^{(R)}v)\bigr)
       +F(P_{\le N}^{(R)}v).
\end{equation*}
For the first term, the pointwise difference bound, weighted H\"older
inequality,
Lemma~\ref{lem:weightedLP}, and Lemma~\ref{lem:weightedSobolev} yield
\begin{align}
 \left\lvert \left\langle F(v)-F(P_{\le N}^{(R)}v),P_N^{(R)}z\right\rangle\right\rvert
 & \quad\lesssim_p
 \mathcal{M}(v)^{\frac{p-1}{q}}
 \left\lVert r^\beta P_{>N}^{(R)}v\right\rVert_{L^{q/2}}
 \left\lVert r^{1/2}P_N^{(R)}z\right\rVert_{L^\infty}\notag\\
 &\quad\lesssim_p
 N^{-4/q}E(v)^{1/2}
 \mathcal{M}(v)^{\frac{p-1}{q}}
 \left\lVert r^{1/2}P_N^{(R)}z\right\rVert_{L^\infty}.
 \label{eq:highvbound}
\end{align}
For the second term, insert a slightly wider projection
$\widetilde P_N^{(R)}$ on $z$ and use the factorization
\begin{equation*}
P_N=N^{-1}\widetilde P_N\sqrt{-\Delta_{\Omega_R}}.
\end{equation*}
The self-adjointness of $\sqrt{-\Delta_{\Omega_R}}$ and
Lemma~\ref{lem:spectralchainHardy} then give
\begin{align}
 \left\lvert \left\langle F(P_{\le N}^{(R)}v),P_N^{(R)}z\right\rangle\right\rvert
 \lesssim_p
 N^{-4/q}E(v)^{1/2}
 \mathcal{M}(v)^{\frac{p-1}{q}}
 \left\lVert r^{1/2}P_N^{(R)}z\right\rVert_{L^\infty}.
 \label{eq:lowvbound}
\end{align}
Indeed, $N^{-1}N^\sigma=N^{-4/q}$. By summing
\eqref{eq:highvbound} and \eqref{eq:lowvbound} over $N$, we obtain
\eqref{eq:paraproduct}, and  completes the proof.
\end{proof}

For the small component $w$, write 
\begin{equation*}
 A(t):=
       \left\lVert w(t)\right\rVert_{L^{3/(1-s_c)}}^2 + \left\lVert w(t)\right\rVert_{L^{3(p-1)}}^{p-1}.
\end{equation*}

\begin{lemma} \label{lem:Taylorremainder}
For every finite-energy radial $v$ and every $w$, we have
\begin{align}
 &\left\lVert \left\lvert v\right\rvert^{p-2}\left\lvert w\right\rvert^2+\left\lvert v\right\rvert\left\lvert w\right\rvert^{p-1}\right\rVert_{L^2}
   \lesssim_p
 E(v)^{1/2}
 \left(\left\lVert w\right\rVert_{L^{a_p}}^2
 +\left\lVert w\right\rVert_{L^{3(p-1)}}^{p-1}\right),
 \label{eq:TaylorL2}
\end{align}
where $ a_p:=\frac3{1-s_c}$. 
Consequently, for $Z\in\{\partial_tv,\nabla v,v/r\}$, we have
\begin{equation}
 \int_{\Omega_R}\left\lvert Z\right\rvert
 \left(\left\lvert v\right\rvert^{p-2}\left\lvert w\right\rvert^2+\left\lvert v\right\rvert\left\lvert w\right\rvert^{p-1}\right)\,\mathrm{d} x
 \lesssim_p E(v)A(t).
 \label{eq:Taylorenergy}
\end{equation}
\end{lemma}

\begin{proof}Let
	\begin{equation*}
	b_p:=\frac{6(p-1)(p-2)}{5p-13}.
	\end{equation*}
The exponent $b_p$ lies between $p+1$ and $6$, and
\begin{equation*}
 \frac12=\frac2{a_p}+\frac{p-2}{b_p}.
\end{equation*}
If
$\theta=(p-3)(p+1)/((p-1)(p-2))$, then
$b_p^{-1}=\theta(p+1)^{-1}+(1-\theta)/6$ and
\begin{equation*}
 \left\lVert v\right\rVert_{L^{b_p}}^{p-2}
 \le \left\lVert v\right\rVert_{L^{p+1}}^{\theta(p-2)}
      \left\lVert v\right\rVert_{L^6}^{(1-\theta)(p-2)}
 \lesssim_p E(v)^{1/2}.
\end{equation*}
Thus H\"older's inequality controls the first term in
\eqref{eq:TaylorL2}. For the
second use
$1/2=1/6+(p-1)/(3(p-1))$ and \eqref{eq:Sobolev6}.
Finally, $\left\lVert Z\right\rVert_2\lesssim E(v)^{1/2}$ for the three displayed choices,
the last one by \eqref{eq:Hardy}. This proves
\eqref{eq:Taylorenergy}.
\end{proof}

\begin{proposition}\label{prop:modifiedenergy}
Let $w$ satisfy \eqref{eq:smallBesov}, and let $v$ be a smooth solution of
\eqref{eq:vequation}. Let
\begin{equation*}
 B(t):=\sum_NN^{-4/q}
 \left\lVert r^{1/2}P_N^{(R)}(\partial_t,\sqrt{-\Delta_{\Omega_R}})w(t)\right\rVert_{L_r^\infty}.
\end{equation*}
Then we have
\begin{equation}
 \frac{\,\mathrm{d}}{\,\mathrm{d} t}\mathcal{E}(t)+c_0\mathcal{M}(v)(t)
 \lesssim_p
 A(t)\mathcal{E}(t)
 +
 \mathcal{E}(t)^{q/4}B(t)^{q/2},
 \label{eq:modifiedinequality}
\end{equation}
where $c_0>0$ depends only on $p$.
\end{proposition}

\begin{proof}
Write
\begin{equation*}
 \mathcal R_2(v,w):=\mathcal R(v,w)-p|v|^{p-1}w.
\end{equation*}
The derivative of the correction in \eqref{eq:modifiedenergy} is
\begin{equation*}
 \frac{\,\mathrm{d}}{\,\mathrm{d} t}\int_{\Omega_R}F(v)w\,\mathrm{d} x
 =p\int_{\Omega_R}|v|^{p-1}\partial_tv\,w\,\mathrm{d} x
 +\int_{\Omega_R}F(v)\partial_tw\,\mathrm{d} x.
\end{equation*}
Combining this identity with \eqref{eq:energyvariation} cancels the
linear-in-$w$ part of the time interaction and gives
\begin{equation*}
 \frac{\,\mathrm{d}}{\,\mathrm{d} t}\left(E(v)+\int_{\Omega_R}F(v)w\,\mathrm{d} x\right)
 =-\int_{\Omega_R}\partial_tv\,\mathcal R_2(v,w)\,\mathrm{d} x
 +\int_{\Omega_R}F(v)\partial_tw\,\mathrm{d} x.
\end{equation*}

We next expand the perturbative term in the Morawetz identity. Radial
integration by parts, using $F(v(t,R))=0$, yields
\begin{align*}
 -cp\int_{\Omega_R}|v|^{p-1}w
 \left(\partial_rv+\frac vr\right)\,\mathrm{d} x& \quad =-c\int_{\Omega_R}w\,\partial_rF(v)\,\mathrm{d} x
 -cp\int_{\Omega_R}F(v)\frac wr\,\mathrm{d} x\notag\\
 &\quad=c\int_{\Omega_R}F(v)\partial_rw\,\mathrm{d} x
 +c(2-p)\int_{\Omega_R}F(v)\frac wr\,\mathrm{d} x.
\end{align*}
Indeed, the coefficient $2$ in the last line is the derivative of the
radial Jacobian $r^2$. Therefore the differentiated identity is
\begin{align}
 \frac{\,\mathrm{d}}{\,\mathrm{d} t}\mathcal{E}(v,w)
 &+\frac{c(p-1)}{p+1}\mathcal{M}(v)
 +2\pi cR^2|\partial_rv(t,R)|^2\notag\\
 &=\mathcal L(t)
 -\int_{\Omega_R}\partial_tv\,\mathcal R_2(v,w)\,\mathrm{d} x -c\int_{\Omega_R}\mathcal R_2(v,w)
 \left(\partial_rv+\frac vr\right)\,\mathrm{d} x,
 \label{eq:modifiedidentity}
\end{align}
where the remaining leading terms have the form
\begin{equation}
 \mathcal L(t):=
 \left\langle F(v),\partial_tw\right\rangle
 +c\left\langle F(v),\partial_rw\right\rangle
 +c(2-p)\left\langle F(v),w/r\right\rangle.
 \label{eq:leadingterms}
\end{equation}
Lemma~\ref{lem:spectralchainHardy} gives
\begin{equation*}
 B_{\partial_rw}(t)+B_{w/r}(t)\lesssim B_{\sqrt{-\Delta_{\Omega_R}} w}(t).
\end{equation*}
Lemma~\ref{lem:paraproduct}, applied to
$z=\partial_tw$, $z=\partial_rw$, and $z=w/r$, therefore gives
\begin{equation*}
 \left\lvert \mathcal L(t)\right\rvert
 \lesssim_p
 E(v)^{1/2}\mathcal{M}(v)^{(p-1)/q}B(t).
\end{equation*}
Young's inequality with conjugate exponents $q/(p-1)$ and $q/2$ yields
\begin{equation}
 E^{1/2}\mathcal{M}^{(p-1)/q}B
 \le \delta\mathcal{M}+C_{p,\delta}E^{q/4}B^{q/2}.
 \label{eq:Youngleading}
\end{equation}
because
\begin{equation*}
 \left(\mathcal{M}^{(p-1)/q}\right)^{q/(p-1)}=\mathcal{M},
 \qquad
 \left(E^{1/2}B\right)^{q/2}=E^{q/4}B^{q/2}.
\end{equation*}

It remains to estimate the Taylor remainders.
Lemma~\ref{lem:Taylorremainder} and \eqref{eq:Taylor2} give
\begin{align*}
 &\int_{\Omega_R}\left\lvert \partial_tv\right\rvert
 \left(\left\lvert v\right\rvert^{p-2}\left\lvert w\right\rvert^2+\left\lvert v\right\rvert\left\lvert w\right\rvert^{p-1}\right)\,\mathrm{d} x
 \lesssim_p E(v)A(t),
 \\
 &\int_{\Omega_R}\left\lvert \nabla v\right\rvert
 \left(\left\lvert v\right\rvert^{p-2}\left\lvert w\right\rvert^2+\left\lvert v\right\rvert\left\lvert w\right\rvert^{p-1}\right)\,\mathrm{d} x
 \lesssim_p E(v)A(t).
\end{align*}
After the leading term has been included in \eqref{eq:leadingterms}, the
$v/r$ remainder satisfies
\begin{equation*}
 \left|\int_{\Omega_R}
 \left(\mathcal R(v,w)-p\left\lvert v\right\rvert^{p-1}w\right)\frac vr\,\mathrm{d} x\right|
 \lesssim_p E(v)A(t).
\end{equation*}
Choose the constant $c$ in \eqref{eq:modifiedenergy} sufficiently small
to ensure the coercivity. Then choose the parameter $\delta$ in
\eqref{eq:Youngleading} so that
\begin{equation*}
 C_p\delta\leq\frac{c(p-1)}{2(p+1)}.
\end{equation*}
The first term on the right of \eqref{eq:Youngleading} is absorbed into
the Morawetz term of \eqref{eq:modifiedidentity}; the boundary flux
is nonnegative and may be discarded. Finally Lemma~\ref{lem:coercivity}
replaces $E$ by $\mathcal{E}$, proving \eqref{eq:modifiedinequality} with
$c_0=c(p-1)/(2(p+1))$.
\end{proof}

By Proposition~\ref{prop:local} and Proposition~\ref{prop:smallBesov},
\begin{equation}
 \int_\mathbb{R} A(t)\,\mathrm{d} t
 \lesssim_p\varepsilon^{p-1}+\varepsilon^2,
 \qquad
 \left\lVert B\right\rVert_{L_t^{q/2}(\mathbb{R})}\lesssim_p\varepsilon.
 \label{eq:ABsmall}
\end{equation}
Indeed, the first summand of $A$ is the $(p-1)$st power of the
$L_t^{p-1}L_x^{3(p-1)}$ norm, the second is the square of the
$L_t^2L_x^{3/(1-s_c)}$ norm, and Minkowski's inequality gives
\begin{equation*}
 \left\lVert B\right\rVert_{L_t^{q/2}}
 \leq\sum_NN^{-4/q}
 \left\lVert r^{1/2}P_N(\partial_t,\sqrt{-\Delta_{\Omega_R}})w\right\rVert_{L_t^{q/2}L_r^\infty}
 =\mathcal{B}_p(w).
\end{equation*}

\section{Proof of main theorem}\label{sec:mainproof}

In this section we prove Theorem~\ref{thm:main}.

\subsection{Normalization of the high--low frequency decomposition}

Let
\begin{equation*}
K:=\left\lVert (u_0,u_1)\right\rVert_{
	\dot H_D^{s_c}(\Omega_1)\times\dot H_D^{s_c-1}(\Omega_1)}.
\end{equation*}
Choose a dyadic $L\gg1$ such that
\begin{equation}
\left\lVert P_{>L}^{(1)}(u_0,u_1)\right\rVert_{
	\dot H_D^{s_c}\times\dot H_D^{s_c-1}}\le\varepsilon,
\label{eq:tailchoice}
\end{equation}
where $\varepsilon$ will be fixed below. Define the rescaled data on
$\Omega_L$ by
\begin{equation*}
U_0(x):=L^{-\frac{2}{p-1}}u_0(x/L),
\qquad
U_1(x):=L^{-\frac{p+1}{p-1}}u_1(x/L).
\end{equation*}
%The critical norm is unchanged and the cut-off $L$ becomes the unit
%cut-off for $-\Delta_{\Omega_L}$.

\begin{lemma}\label{lem:scalingcutoff}
For $L\ge1$, let 
\begin{equation*}
 (\mathcal S_Lf)(x):=L^{-2/(p-1)}f(x/L),\qquad x\in\Omega_L.
\end{equation*}
Then we have
\begin{equation}
 \sqrt{-\Delta_{\Omega_L}}\,\mathcal S_L
 =L^{-1}\mathcal S_L\sqrt{-\Delta_{\Omega_1}},
 \qquad
 P_{\le1}^{(L)}\mathcal S_L
 =\mathcal S_LP_{\le L}^{(1)}.
 \label{eq:cutoffscaling}
\end{equation}
The $\dot H_D^{s_c}$ norm is invariant under $\mathcal S_L$. If choosing 
$$V_0=P_{\le1}^{(L)}\mathcal S_Lu_0, \;
V_1=P_{\le1}^{(L)}L^{-1}\mathcal S_Lu_1, $$ then we have
\begin{equation}
 \left\lVert \nabla V_0\right\rVert_2+\left\lVert V_1\right\rVert_2\lesssim K,
 \qquad
 \left\lVert V_0\right\rVert_{p+1}^{p+1}\lesssim_p K^{p+1},
 \label{eq:cutoffenergybounds}
\end{equation}
with constants independent of $L$.
\end{lemma}

\begin{proof}
The first identity follows either by differentiation or from
\eqref{eq:conjugation}. In deed,  by a change of variables, we have 
\begin{equation*}
 -\Delta_{\Omega_L}\mathcal S_Lf
 =L^{-2}\mathcal S_L(-\Delta_{\Omega_1}f).
\end{equation*}
The spectral theorem then yields both identities in
\eqref{eq:cutoffscaling}. 

For every $s\in\mathbb{R}$, we have
\begin{equation}
 \left\lVert \mathcal S_Lf\right\rVert_{\dot H_D^s(\Omega_L)}
 =L^{\frac32-s-\frac{2}{p-1}}
 \left\lVert f\right\rVert_{\dot H_D^s(\Omega_1)}.
 \label{eq:Dirichletnormscaling}
\end{equation}
For $s=s_c$, the exponent vanishes. Similarly,
\begin{equation*}
 \left\lVert L^{-1}\mathcal S_Lg\right\rVert_{\dot H_D^{s_c-1}(\Omega_L)}
 =\left\lVert g\right\rVert_{\dot H_D^{s_c-1}(\Omega_1)}.
\end{equation*}

Because $\lambda\leq2$ on the frequency support of $P_{\leq1}^{(L)}$, we have
\begin{align*}
 \left\lVert \nabla V_0\right\rVert_2^2
 &=4\pi\int_0^2\lambda^2|\mathcal F_LV_0(\lambda)|^2\,\mathrm{d}\lambda
 \lesssim\left\lVert V_0\right\rVert_{\dot H_D^{s_c}}^2,\\
 \left\lVert V_1\right\rVert_2^2
 &=4\pi\int_0^2|\mathcal F_LV_1(\lambda)|^2\,\mathrm{d}\lambda
 \lesssim\left\lVert V_1\right\rVert_{\dot H_D^{s_c-1}}^2.
\end{align*}
This proves the first estimate in \eqref{eq:cutoffenergybounds}, namely,
that the low-frequency data have finite kinetic energy. For the potential
term, the Littlewood--Paley projection $P_{\leq1}^{(L)}$ is bounded on
$L^{p+1}$ by the Bernstein estimate in Lemma $2.3$ in \cite{XuYang2024}, and hence we have
\begin{equation*}
 \left\lVert V_0\right\rVert_{p+1}^{p+1}
 \leq  \left\lVert \mathcal S_Lu_0\right\rVert_{p+1}^{p+1}
 =L^{3-\frac{2(p+1)}{p-1}}\left\lVert u_0\right\rVert_{p+1}^{p+1}
 =L^{\frac{p-5}{p-1}}\left\lVert u_0\right\rVert_{p+1}^{p+1}.
\end{equation*}
Because $L\geq1$ and $p<5$, the right-hand side is no larger than
$\left\lVert u_0\right\rVert_{p+1}^{p+1}$. This implies the second estimate in \eqref{eq:cutoffenergybounds} by the Sobolev embedding inequality \eqref{eq:criticalLp} on
$\Omega_1$.
\end{proof}

By \eqref{eq:cutoffscaling}, the rescaled decomposition is
\begin{align*}
 (V_0,V_1)
 &=P_{\leq1}^{(L)}(U_0,U_1)
 =\left(\mathcal S_LP_{\leq L}^{(1)}u_0,
 L^{-1}\mathcal S_LP_{\leq L}^{(1)}u_1\right),\\
 (W_0,W_1)
 &=P_{>1}^{(L)}(U_0,U_1)
 =\left(\mathcal S_LP_{>L}^{(1)}u_0,
 L^{-1}\mathcal S_LP_{>L}^{(1)}u_1\right).
\end{align*}
Thus \eqref{eq:tailchoice} and
\eqref{eq:Dirichletnormscaling} imply
\begin{equation*}
 \left\lVert (W_0,W_1)\right\rVert_{\dot H_D^{s_c}(\Omega_L)
 \times\dot H_D^{s_c-1}(\Omega_L)}\leq\varepsilon,
\end{equation*}
while the low-frequency component $(V_0,V_1)$  has finite energy. More precisely, we have 
\begin{equation*}
 E(V)(0)\lesssim_p K^2+K^{p+1}=:E_K.
\end{equation*}
This estimate is independent of $L$. 

\subsection{Global analysis of the difference propagation with the low-frequency part}

Let $W$ be the small nonlinear solution associated with the high-frequency
data, and let $V$ solve the difference equation.
Set
\begin{equation*}
 \alpha:=\frac q4=\frac{p+1}{4}\in \left(1, \frac32\right).
\end{equation*}
Let $I_+=[0,T_+)$ be the maximal forward lifespan of the  energy
solution $V$, and define
\begin{equation*}
 \mathfrak T:=\left\{T\in I_+:
 \sup_{0\leq t\leq T}\mathcal{E}(t)\leq4C_pE_K\right\}.
\end{equation*}
Coercivity and continuity show that $\mathfrak T$ contains a
neighborhood of zero after increasing $C_p$ if necessary. On a connected
component $[0,T]\subset\mathfrak T$,  we drop the positive Morawetz term in
\eqref{eq:modifiedinequality} and apply Lemma~\ref{lem:Bihari}, with
\begin{equation*}
 a(t)=C_pA(t),\qquad b(t)=C_pB(t)^{2\alpha}.
\end{equation*}
Since \eqref{eq:ABsmall} implies
\begin{equation*}
 \exp\left((\alpha-1)\int_0^Ta(t)\,\mathrm{d} t\right)
 \leq\exp\bigl(C_p(\varepsilon^{p-1}+\varepsilon^2)\bigr)
 \leq2
\end{equation*}
for small $\varepsilon$, Lemma~\ref{lem:Bihari} gives
\begin{equation}
 \mathcal{E}(t)^{1-\alpha}
 \geq\frac12\mathcal{E}(0)^{1-\alpha}
 -C_p\int_0^tB(s)^{2\alpha}\,\mathrm{d} s.
 \label{eq:Bihari}
\end{equation}
By \eqref{eq:ABsmall}, we have
\begin{equation*}
 \int_\mathbb{R} B(t)^{2\alpha}\,\mathrm{d} t
 \lesssim_p\varepsilon^{2\alpha}.
\end{equation*}
If choosing the high frequency truncation in \eqref{eq:tailchoice} so small that
\begin{equation}
 C_p(1+E_K)^{\alpha-1}\varepsilon^{2\alpha}\le\frac18.
 \label{eq:epsilonchoice}
\end{equation}
By \eqref{eq:Bihari}--\eqref{eq:epsilonchoice}, we obtain the improved estimate
\begin{equation*}
 \sup_{0\leq t\leq T}\mathcal{E}(t)\leq2C_pE_K
\end{equation*}
since $1-\alpha<0$. 
Thus $\mathfrak T$ is both open and closed in $I_+$. The continuation
criterion in Proposition~\ref{prop:affine} rules out $T_+<\infty$.
Repeating the argument for the negative time direction gives
\begin{equation}
 \sup_{t\in\mathbb{R}}E(V)(t)\lesssim_{p,K}1.
 \label{eq:globalenergy}
\end{equation}
Thus the energy solution $V$ exists globally on $\Omega_L$.

\subsection{Persistence of critical regularity}

\begin{lemma}\label{lem:persistence}
Let
\begin{equation*}
 \delta_c:=1-s_c,
 \qquad
 \delta_*:=\frac{3(p-1)}{2(p+1)}.
\end{equation*}
Then $0<\delta_c<\delta_*<1$. If $V$ is an energy solution,
$W\in C_t\dot H_D^{s_c}$, and
$(V_0,V_1)\in\dot H_D^{s_c}\times\dot H_D^{s_c-1}$, then we have 
\begin{equation*}
 (V,\partial_tV)\in
 C_{\mathrm{loc}}(\dot H_D^{s_c}\times\dot H_D^{s_c-1}).
\end{equation*}
\end{lemma}

\begin{proof}
Set $\mathcal N(V,W)=F(V+W)-F(W)=F(V)+\mathcal R(V,W)$. Since
$6/(3+2\delta_*)=(p+1)/p$, Lemma~\ref{lem:basictools} and H\"older's
inequality imply
\begin{align}
 \left\lVert \mathcal N(V,W)\right\rVert_{\dot H_D^{-\delta_*}}
 &\lesssim_p
 \left\lVert \mathcal N(V,W)\right\rVert_{L^{(p+1)/p}} \lesssim_p
 \left\lVert V\right\rVert_{L^{p+1}}^p
 +\left\lVert V\right\rVert_{L^{p+1}}\left\lVert W\right\rVert_{L^{p+1}}^{p-1}.
 \label{eq:nonlinearnegative}
\end{align}
The right-hand side is locally bounded by the energy of $V$ and the
critical Sobolev norm of $W$.

Set $H(t)=V(t)-V_0$. The energy solution satisfies
$H\in C_{\mathrm{loc}}L^2\cap C_{\mathrm{loc}}\dot H_D^1$. Indeed,
\begin{equation*}
 H(t)-H(t_0)=\int_{t_0}^t\partial_sV(s)\,\mathrm{d} s,
\end{equation*}
so the $L^2$ continuity follows from
$\partial_tV\in C_{\mathrm{loc}}L^2$, while the $\dot H_D^1$
continuity is part of the construction in the energy space. Applying
\eqref{eq:spectralinterpolation} to the difference
$H(t)-H(t_0)$ gives
\begin{equation*}
 \left\lVert H(t)-H(t_0)\right\rVert_{\dot H_D^{s_c}}
 \leq\left\lVert H(t)-H(t_0)\right\rVert_2^{1-s_c}
 \left\lVert H(t)-H(t_0)\right\rVert_{\dot H_D^1}^{s_c},
\end{equation*}
and hence $H\in C_{\mathrm{loc}}\dot H_D^{s_c}$. For the estimate of $ \partial_tV$, use
\begin{align}
 \partial_tV(t)
 &=\cos(t\sqrt{-\Delta_{\Omega_R}})V_1
 -\sqrt{-\Delta_{\Omega_R}}\sin(t\sqrt{-\Delta_{\Omega_R}})V_0\notag\\
 &\qquad-\int_0^t\cos((t-s)\sqrt{-\Delta_{\Omega_R}})
 \mathcal N(V,W)(s)\,\mathrm{d} s.
 \label{eq:velocityDuhamel}
\end{align}

On $\{\lambda\le1\}$, the embedding
$\dot H_D^{-\delta_*}\hookrightarrow\dot H_D^{-\delta_c}$ follows
directly from \eqref{eq:spectralnorm}:
\begin{equation*}
 \left\lVert P_{\leq1}f\right\rVert_{\dot H_D^{-\delta_c}}^2
 =4\pi\int_0^1\lambda^{-2\delta_c}|\mathcal F_Rf|^2\,\mathrm{d}\lambda
 \leq\left\lVert P_{\leq1}f\right\rVert_{\dot H_D^{-\delta_*}}^2.
\end{equation*}
The local boundedness in \eqref{eq:nonlinearnegative} and dominated
convergence in the time integral show that the low-frequency part of
\eqref{eq:velocityDuhamel} is continuous in
$\dot H_D^{s_c-1}=\dot H_D^{-\delta_c}$. 

On $\{\lambda>1\}$, we have
\begin{equation*}
 \left\lVert P_{>1}f\right\rVert_{\dot H_D^{s_c-1}}
 \leq\left\lVert P_{>1}f\right\rVert_2,
\end{equation*}
because $s_c-1<0$. The energy continuity
$\partial_tV\in C_{\mathrm{loc}}L^2$ therefore supplies the
high-frequency continuity. Combining the two regions proves the result.
\end{proof}

On every compact interval $[-T,T]$, write
\begin{equation*}
 V(t)=V_0+H(t),
 \qquad
 H(t)=\int_0^t\partial_sV(s)\,\mathrm{d} s.
\end{equation*}
The energy bound gives
\begin{equation*}
 \sup_{\left\lvert t\right\rvert\le T}\left\lVert H(t)\right\rVert_2
 \lesssim_{p,K}T,
 \qquad
 \sup_{\left\lvert t\right\rvert\le T}\left\lVert \nabla H(t)\right\rVert_2
 \lesssim_{p,K,T}1.
\end{equation*}
Lemma~\ref{lem:persistence} now gives critical continuity for $V$.
Since $V_0,V_1$ already have the critical regularity, and $W$ satisfies
\eqref{eq:smalldata},
\begin{equation*}
 (U,\partial_tU)=(V+W,\partial_tV+\partial_tW)
 \in C\!\left(\mathbb{R};\dot H_D^{s_c}\times\dot H_D^{s_c-1}\right).
\end{equation*}

It remains to justify the passage from smooth solutions to critical data.
Choose smooth spectral
truncations
\begin{equation*}
 (U_{0,n},U_{1,n})
 :=\bigl(P_{[1/n,n]}^{(L)}U_0,P_{[1/n,n]}^{(L)}U_1\bigr),
\end{equation*}
where $P_{[a,b]}^{(L)}$ denotes a smooth spectral multiplier equal to one
on $[2a,b/2]$ and supported in $[a,b]$. Follow this, if necessary, by a
compactly supported radial approximation in
the sine-transform variable. Then
\begin{equation*}
 (U_{0,n},U_{1,n})\longrightarrow(U_0,U_1)
 \quad\hbox{in }\dot H_D^{s_c}\times\dot H_D^{s_c-1}.
\end{equation*}
For all sufficiently large $n$, the same unit cut-off makes the
high-frequency tails smaller than $2\varepsilon$. The estimates above
are consequently uniform in $n$. On a fixed interval $[-T,T]$, divide
time into finitely many subintervals on which the critical spacetime norm
of one approximate solution is below the stability threshold in
Corollary~\ref{cor:stability}. Iteration of
\eqref{eq:stabilityestimate} shows that
\begin{equation*}
 (U_n,\partial_tU_n)\longrightarrow(U,\partial_tU)
 \quad\hbox{in }C([-T,T];
 \dot H_D^{s_c}\times\dot H_D^{s_c-1}).
\end{equation*}
The limit satisfies the Duhamel formula, is independent of the
approximating sequence, and is therefore the unique critical mild
solution.

Finally consider two initial data sets in a sufficiently small
neighborhood of a fixed datum. By \eqref{eq:spectralTail}, the same $L$
in \eqref{eq:tailchoice} makes both high-frequency tails small. The
difference estimates in Lemma~\ref{lem:criticalnonlinear}
and Lemma~\ref{lem:energyperturbation}, followed by the same finite subdivision of
$[-T,T]$, give
\begin{align*}
 &\left\lVert (U-\widetilde U,\partial_tU-\partial_t\widetilde U)\right\rVert_
 {C([-T,T];\dot H_D^{s_c}\times\dot H_D^{s_c-1})} \leq C_{p,K,T}
 \left\lVert (U_0-\widetilde U_0,U_1-\widetilde U_1)\right\rVert_
 {\dot H_D^{s_c}\times\dot H_D^{s_c-1}}.
\end{align*}
This proves local Lipschitz dependence. Scale the solution on
$\Omega_L$ back by
\begin{equation*}
 u(t,x)=L^{\frac{2}{p-1}}U(Lt,Lx).
\end{equation*}
The resulting solution is defined on $\Omega_1$ and satisfies
\eqref{eq:mainconclusion}. This completes the proof of
Theorem~\ref{thm:main}.

%  \noindent \subsection*{Acknowledgements.}
%G.~Xu gratefully acknowledges support from NSFC (No. 12371240, No. 12431008).  

% The author would like to thank the referee for his/her valuable comments and suggestions, which helps to improve our manuscript.
%


\begin{thebibliography}{99}
	
	\bibitem{BahChe:NLW}
	H.~Bahouri, J.~Y.~Chemin,
	\newblock{On global well-posedness for defocusing cubic wave equation}, 
	\newblock{\em Int. Math. Res. Not. IMRN},   2006, Art. ID 54873, 12 pp.

\bibitem{Bourgain1998}
J.~Bourgain,
\newblock{Refinements of Strichartz' inequality and applications to 2D-NLS
with critical nonlinearity},
\newblock{\em Int. Math. Res. Not. IMRN}, \textbf{1998} (1998), 253--283.

\bibitem{Bourgain1999}
J.~Bourgain,
\newblock{Global wellposedness of defocusing critical nonlinear Schr\"odinger equation in the radial case}. 
\newblock{\em 
J. Amer. Math. Soc.}, \textbf{12} (1999), no. 1, 145–-171.

  \bibitem{Burq03: NLW CPDE}
N.~Burq,
\newblock{Global Strichartz estimates for nontrapping geometries: about an article by H. Smith and C. Sogge}.
\newblock{\em Comm. Partial Differential Equations}, \textbf{28}(2003), no. 9-10, 1675--1683.


\bibitem{BurqLP:NLW}
N.~Burq, G.~Lebeau, F.~Planchon,
\newblock{Global existence for energy critical waves in 3-D domains}. 
\newblock{\em J. Amer. Math. Soc.}, \textbf{ 21}(2008), no. 3., 831--845.

\bibitem{BurqP:def NLW: Neum prob}
N.~Burq, F.~Planchon,
\newblock{Global existence for energy critical waves in 3-D domains: Neumann boundary conditions}. 
\newblock{\em Amer. J. Math.},  \textbf{131}(2007), 1715--1742.

\bibitem{CKSTT2008}
J.~Colliander, M.~ Keel, G.~ Staffilani, H.~ Takaoka, and T.~ Tao, 
\newblock{Global well-posedness and scattering for the energy-critical nonlinear Schrödinger equation in  $\mathbb{R}^3$}.
\newblock{\em
Ann. of Math.}, \textbf{167} (2008), no. 3, 767–865.

\bibitem{DodsonCubic2019}
B.~Dodson,
\newblock{Global well-posedness for the defocusing, cubic, nonlinear wave equation in three dimensions for radial initial data in $\dot H^{s}\times \dot H^{s-1}$,  $s>1/2$}.
\newblock{\em Int. Math. Res. Not. IMRN}, \textbf{2019}(2019), no. 21, 6797-–6817.

\bibitem{DodsonCubic2018}
B.~Dodson,
\newblock{Global well-posedness and scattering for the radial, defocusing,
	cubic nonlinear wave equation},
\newblock{\em Duke Math. J.}, \textbf{170}(2021), no. 15, 3267--3321. 

\bibitem{Dodson2018}
B.~Dodson,
\newblock{Global well-posedness for the radial, defocusing, nonlinear wave
equation for $3<p<5$},
\newblock{\em Amer. J. Math.}, \textbf{146} (2024), no. 1, 1--46.

\bibitem{Dodson2024}
B.~Dodson,
\newblock{Sharp global well-posedness and scattering for the radial conformal nonlinear wave equation}.
\newblock{\em Int. Math. Res. Not. IMRN}, \textbf{2024}(2024), no. 9, 7748–-7775.

\bibitem{Dodson2025}
B.~Dodson,
\newblock{Sharp global well-posedness and scattering for the radial nonlinear intercritical wave equation}. 
\newblock{\em Indiana Univ. Math. J.} \textbf{74} (2025), no. 6, 1677–-1720.

\bibitem{DuyLaf:NLW:Neum prob}
T.~Duyckaerts, D.~Lafontaine,
\newblock{Scattering for critical radial Neumann waves outside a ball}.
\newblock{\em Rev. Mat. Iberoam.}, \textbf{38}(2022), no. 2, 659--703. 

\bibitem{DuyYang24:LW Str}
T.~Duyckaerts, J. Yang,
\newblock{Dispersive estimates for wave and Schr\"odinger equations with a potential in non-trapping exterior domains}.
\newblock{\em arXiv:2401.12608}, 2024.


\bibitem{GallPlan03:NLW}
I.~Gallagher, F.~ Planchon,
\newblock{On global solutions to a defocusing semi-linear wave equation}. 
\newblock{\em Rev. Mat. Iberoamericana}, \textbf{19}(2003),  no. 1, 161--177.

\bibitem{GinV:NLW:St est}
J.~Ginibre, G.~Velo,
\newblock{Generalized Strichartz inequalities for the wave equation}.
\newblock{\em  J. Funct. Anal.}, \textbf{133}(1995), 50--68.

\bibitem{IvLebeauP:LW:dis est I}
O.~Ivanovici, G.~Lebeau, F.~Planchon,
\newblock{Dispersion for the wave equation inside strictly convex domains I: the Friedlander model case}.
\newblock{\em Ann. Math.}, \textbf{180}(2014), 323--380.

\bibitem{IvLasLebP:LW:dis est II}
O.~Ivanovici, R.~Lascar, G.~Lebeau, F.~Planchon,
\newblock{Dispersion for the wave equation inside strictly convex domains II: the general case}.
\newblock{\em Ann. PDE.}, \textbf{9}(2023), no. 2, Paper No. 14, 117 pp.


\bibitem{KeelTao1998}
M.~Keel and T.~Tao,
\newblock{Endpoint Strichartz estimates},
\newblock{\em Amer. J. Math.}, \textbf{120} (1998), 955--980.

\bibitem{KenigPV00:NLW}
C.~E.~Kenig, G.~Ponce, L.~Vega,
\newblock{Global well-posedness for semi-linear wave equations}.
\newblock{\em  Comm. Partial Differential Equations}, 25(2000), no. 9-10, 1741--1752.

\bibitem{KenigM2006:NLS}
C.~E.~Kenig, F.~Merle, 
\newblock{
Global well-posedness, scattering and blow-up for the energy-critical, focusing, non-linear Schr\"dinger equation in the radial case}.
\newblock{\em Invent. Math.}, \textbf{166} (2006), no. 3, 645–-675.

\bibitem{KenigM2008:NLW}
C.~E.~Kenig, F.~Merle, 
\newblock{Global well-posedness, scattering and blow-up for the energy-critical focusing non-linear wave equation}.
\newblock{\em Acta Math.}, \textbf{201} (2008), no. 2, 147–-212.

\bibitem{KillpV2008}
R.~Killip and M.~Vsian,
\newblock{Nonlinear Schr\"odinger equations at critical regularity}.
\newblock{Clay Math. Proceedings}, \textbf{ 17},
Amer. Math. Soc., Providence, RI, 2013. 

\bibitem{KVZ16:IMRN Riesz tf}
R.~Killip, M.~Visan, X.~Zhang,
\newblock{Riesz transforms outside a convex obstacle}.
\newblock{\em Int. Math. Res. Not. IMRN},  Vol. \textbf{2016}(2016), no. 19, pp. 5875--5921. 

\bibitem{LiSZ:NLS}
D.~Li, H.~Smith, X.~Zhang,
\newblock{Global well-posedness and scattering for defocusing energy-critical NLS in the exterior of balls with radial data}.
\newblock{\em Math. Res. Lett.}, \textbf{19}(2012), no. 1, 213--232.

\bibitem{LiXZ:NLS}
D.~Li, G.~Xu, X.~Zhang,
\newblock{On the dispersive estimate for the Dirichlet Schr\"odinger propagator and applications to energy critical NLS}.
\newblock{\em Canad. J. Math.}, \textbf{66}(2014), no. 5, 1110--1142.


\bibitem{Lin2026}
Z.~Lin,
\newblock{Low-regularity global well-posedness for the defocusing cubic wave
equation outside a ball with radial initial data},
arXiv:2610.38495, 2026.

 \bibitem{Roy:NLW:radial}
T.~Roy, 
\newblock{Global well-posedness for the radial defocusing cubic wave equation on $\R^3$ and for rough data}, 
\newblock{\em Electron. J. Differential Equations},   2007, No. 166, 22 pp.

\bibitem{Roy:NLW}
T.~Roy,
\newblock{Adapted linear-nonlinear decomposition and global well-posedness for solutions to the defocusing cubic wave equation on $\R^3$}, 
\newblock{\em Discrete Contin. Dyn. Syst.},  24 (2009), no. 4, 1307--1323.

      
\bibitem{SmithSogge95:NLW JAMS}
H.~F.~Smith, C.~D.~Sogge,
\newblock{On the critical semilinear wave equation outside convex obstacles}.
\newblock{\em  J. Amer. Math. Soc.},  8(1995), no. 4, 879--916. 

  \bibitem{SmithSogge00:NLW CPDE}
H.~F.~Smith, C.~D.~Sogge,
\newblock{Global Strichartz estimates for nontrapping perturbations of the Laplacian: Estimates for nontrapping perturbations}.
\newblock{\em Comm. Partial Differential Equations }, 25(2000), no. 11-12, 2171--2183.


\bibitem{SteinWeiss1958}
E.~M.~Stein and G.~Weiss,
\newblock{Fractional integrals on $n$-dimensional Euclidean space},
\newblock{\em J. Math. Mech.}, \textbf{7} (1958), 503--514.


\bibitem{TaoV2005}
T.~Tao and M.~Visan,
\newblock{Stability of energy-critical nonlinear Schrödinger equations in high dimensions}. 
\newblock{\em
	Electron. J. Differential Equations},  2005, No. 118, 28 pp.



\bibitem{XuYang2024}
G.~Xu and P.~Yang,
\newblock{Global well-posedness of the defocusing, cubic nonlinear wave
equation outside of the ball with radial data},
\newblock{\em Trans. Amer. Math. Soc.}, \textbf{379} (2026), no. 9,  6661--6685.

\bibitem{XuYangYou2025}
G.~Xu, P.~Yang, and Z.~You,
\newblock{Global well-posedness of the defocusing nonlinear wave equation
outside of a ball with radial data for $3<p<5$},
arXiv:2501.05188, 2025.


\end{thebibliography}
\end{document}